\documentclass{article} 
\usepackage{amssymb}
\usepackage{amsmath}
\usepackage{latexsym}
\usepackage[mathscr]{euscript}
\usepackage{graphicx}
\usepackage{amscd}
\usepackage{amsthm} 
\usepackage{fullpage} 
\usepackage{wrapfig}
\newtheorem{theorem}{Theorem}
\newtheorem{corollary}[theorem]{Corollary}
\newtheorem{lemma}[theorem]{Lemma} 
\newtheorem{proposition}[theorem]{Proposition}

\newtheorem{definition}{Definition}
\newcommand{\bc}{\mathbb{C}}

\newcommand{\bz}{\mathbb{Z}}
\newcommand{\br}{\mathbb{R}}

\newcommand{\bh}{\mathbb{H}}

\newcommand{\bs}{\mathbb{S}}

\newcommand{\p}{\partial}
\newcommand{\cc}{\mathcal{C}}
\newcommand{\cd}{\mathcal{D}}
\newcommand{\ce}{\mathcal{E}}

\newcommand{\cl}{\mathcal{L}}

\newcommand{\cg}{\mathcal{G}}
\newcommand{\cb}{\mathcal{B}}
\newcommand{\cp}{\mathcal{P}}

\newcommand{\cm}{\mathcal{M}}

\newcommand{\cf}{\mathcal{F}}

\newcommand{\ca}{\mathcal{A}}

\newcommand{\cj}{\mathcal{J}}
\newcommand{\hk}{\hookrightarrow}
\newcommand{\bg}{\bigskip}
\newcommand{\med}{\medskip}
\newcommand{\la}{\longrightarrow}
\newcommand{\bfl}{\begin{flushleft}}
\newcommand{\efl}{\end{flushleft}}

\newcommand{\eps}{\epsilon}

\newcommand{\colim}{\operatorname{colim}}

\newcommand{\xr}{\xrightarrow}

\newcommand{\bcm}{\bar \cm}

 \newcommand{\ut}{\underbar{2}}

\begin{document}

\title{ Floer homotopy theory, revisited} 
\author{Ralph L. Cohen \thanks{The author was partially supported by a grant from the NSF.} \\ Department of Mathematics \\Stanford University \\ Bldg. 380 \\ Stanford, CA 94305}
\date{\today}
%\subjclass{ 57R50; 30F99; 57M07}
\maketitle 
%\begin{centerline} {\sl Dedicated to the memory of Samuel Gitler. }
%\end{centerline}
\begin{abstract} In 1995 the author, Jones, and Segal introduced the notion of ``Floer homotopy theory" \cite{cjs}. The proposal was to attach a (stable) homotopy type to the geometric data given in a version of Floer homology.  More to the point, the question was asked, ``When is the Floer homology isomorphic to the (singular) homology of a naturally occuring (pro)spectrum  defined from the properties of the moduli spaces inherent in the Floer theory?". A proposal for how to construct such a spectrum was given in terms of a ``framed flow category", and some rather simple examples were described.  Years passed before this notion found some genuine applications to  symplectic geometry and low dimensional topology.  However in recent years several striking applications have been found, and the theory has been developed on a much deeper level. Here we summarize some of these exciting developments,      and describe some of the new techniques that were introduced.  Throughout we try to point out that this area is a very fertile ground at the interface of homotopy theory, symplectic geometry, and low dimensional topology.  \end{abstract}

\tableofcontents

\section*{Introduction} In three seminal papers in 1988 and 1989    A. Floer introduced Morse theoretic homological invariants  that transformed the study of low dimensional topology and  symplectic geometry.  In \cite{floerI} Floer defined an ``instanton homology" theory for 3-manifolds that, when paired with Donaldson's polynomial invariants of 4 manifolds defined a gauge theoretic 4-dimensional topological field theory that revolutionized the study of low dimensional topology and geometry. In \cite{floerFix}, Floer defined an infinite dimensional Morse theoretic homological invariant for symplectic manifolds, now referred to as  ``Symplectic" or ``Hamiltonian"  Floer homology,  that allowed him to prove a well-known conjecture of Arnold on the number of fixed points of a diffeomorphism $\phi_1 : M \xr{\cong} M$
arising from a time-dependent Hamiltonian flow $\{\phi_t\}_{0 \leq t \leq 1}$.  In \cite{floerLag} Floer introduced   ``Lagrangian intersection Floer theory"   for the study of interesections of Lagrangian submanifolds of a symplectic manifold.  

Since that time there have been many other versions of Floer theory introduced in geometric topology, including a 
Seiberg-Witten Floer homology \cite{KM}.  This is similar in spirit to Floer's ``instanton homology", but it is based on the Seiberg-Witten equations rather than the Yang-Mills equations.  There were many difficult, technical analytic issues in developing Seiberg-Witten Floer theory, and Kronheimer and Mrowka's book \cite{KM} deals with them masterfully and elegantly.  Another important geometric theory is 
 Heegard Floer theory introduced by Oszvath and Szabo \cite{OS}.  This is an invariant  of a closed 3-manifold equipped with a $spin^c$ structure. It is computed using a Heegaard diagram of the manifold.  It allowed for  a related ``knot Floer homology" introduced by  Oszvath and Szabo \cite{OS2} and by Rasmussen \cite{Ras}.  Khovanov's important homology theory that gave a ``categorification" of the Jones polynomial \cite{khov}    was eventually   was shown to be related to Floer theory by Seidel and Smith \cite{SS} and Abouzaid and Smith \cite{AbS}.  Lipshitz and Sarkar  \cite{lipsark} showed that there is an associated ``Khovanov stable homotopy".   There have been many other variations of  Floer  theories as well. 

 The rough idea in all of these theories is to associate a Morse-like chain complex generated by the critical points of a functional defined typically on an infinite dimensional space. Recall that in classical Morse theory, given a Morse function $f : M\to \br$ on a closed Riemannian manifold, the ``Morse complex"  is the chain complex
 $$
 \cdots \to C_p(f)  \xr{\p_p} C_{p-1}(f) \to \cdots 
 $$
 where $C_p(f)$ is the free abelian group generated by the critical points of $f$ of index $p$, and the boundary homomorphisms can be computed by ``counting" the gradient flow-lines connecting critical points of relative index one.
 More specfically,  if $Crit_q(f)$ is the set of critical points of $f$ of index $q$, then if $a \in Crit_p(f)$, then  
\begin{equation}\label{morsecomp}
 \p_p ([a])  =  \sum_{b \in Crit_{p-1}(f)}  \#\cm(a,b) \, [b]
\end{equation}
 where $\cm (a,b)$ is the moduli  space of gradient flow lines connecting $a$ to $b$, which, since $a$ and $b$ have relative index one is a closed, zero dimensional oriented manifold.  $\#\cm (a,b)$ reflects the ``oriented count" of this finite set.  More carefully  $\#\cm (a,b)$ is the integer in the zero dimensional  oriented cobordism group, $\pi_0 MSO \cong \bz$ represented by $\cm(a, b)$.      
 
   In Floer's original examples, the functionals he studied were in fact $\br/\bz$-valued.  In the case of Floer's instanton theory, the  relevant functional is the  Chern-Simons map defined on the space of connections on   a principal $SU(2)$-bundle over the three-manifold. Its critical points are flat connections and its flow lines are ``instantons", i.e   anti-self-dual connections on the three-manifold $Y$ crossed with the real line. Modeling classical Morse theory, the ``Floer complex" is generated by the critical points of this functional, suitably perturbed to make them nondegenerate, and the boundary homomorphisms are computed by  taking oriented counts of the gradient flow lines, i.e anti-self-dual connections on $Y \times \br$ that connect critical points of relative index one.    
   
   When Floer introduced what is now called ``Symplectic" or ``Hamiltonian" Floer homology (``SFH") to use in his proof of Arnold's conjecture,  he studied the symplectic action defined on the free loop space of the underlying symplectic manifold $M$,  $$ \ca : LM \to \br/\bz. $$    He perturbed the action functional by a time dependent Hamiltonian function  $H : \br/\bz \times M \to \br$.  Call the resulting functional $\ca_H$.  The critical  points of $\ca_H$ are the 1-periodic orbits of the Hamiltonian vector field.  That is, they are smooth loops $\alpha : \br/\bz \to M $ satisfying the differential equaiton
   $$
   \frac{d\alpha}{dt}  = X_H(t, \alpha (t))
   $$
   where $X_H$ is the Hamiltonian vector field.  One way to think of $X_H$ is that the symplectic 2-form $\omega$ on $M$ defines, since it is nondegenerate, a bundle isomorphism $\omega : TM \xr{\cong} T^*M$.  This   induces  an identification of the periodic one-forms, that is sections of the cotangent bundle pulled back over $\br/\bz \times M$,  with periodic vector fields.
   The Hamiltonian vector field $X_H$ corresponds to the differential $dH$ under this identification.  Floer showed that with respect to a generic Hamiltonian,  the critical points  of $\ca_H$ are nondegenerate.
   
   Of course to understand the gradient flow lines connecting critical points, one must have a metric.  This is defined using the symplectic form and a  compatible choice of almost complex structure $J$ on $TM$.  With respect to this structure the gradient flow lines are curves $\gamma : \br \to LM$, or equivalently, maps of cylinders
   $$
   \gamma : \br \times S^1 \to M
   $$
 that satisfy  (perturbed) Cauchy-Riemann equations.   If $\cm (\alpha_1, \alpha_2,  H, J)$ is the moduli space of such ``pseudo-holomorphic" cylinders that connect the periodic orbits $\alpha_1$ and $\alpha_2$, then the boundary homomorphisms Morse-Floer chain complex is defined by giving an oriented count of the zero dimensional moduli spaces, in analogy to the situation in classical Morse theory described above. 
 
 The other, newer examples of Floer homology tend to be similar.  There is typically a functional that can be perturbed in such a way that its critical points and zero dimensional moduli spaces of gradient flow lines define a chain complex whose homology is invariant of the choices made.
 
 In the case of a  classical Morse  function $f : M \to \br$ closed manifold,  the Morse complex can be viewed as the cellular chain complex of a $CW$-complex $X_f$ of the homotopy type of $M$.  $X_f$ has one cell of dimension $\lambda$, for every critical point of $f$ of index $\lambda$.  The attaching maps were studied by Franks in \cite{franks}, and one may view the work of the author and his collaborators in \cite{cjs1}, \cite{cjs} as a continuation of that study.   This led us to ask the  question:

\medskip
\noindent \bf Q1: \rm  Does the Floer chain complex arise as the cellular chain complex of a $CW$-complex or a $C.W$-spectrum? 

More specifically we asked the the following question: 

\medskip
\noindent \bf Q2: \rm  What properties of the data of a Floer functional, i.e its critical points and moduli spaces of gradient flow lines connecting them,  are needed to define a $CW$-spectrum realizing the Floer chain complex?

More generally one might ask the following question:

\medskip
\noindent \bf Q3:  \rm  Given a finite chain complex,  is there a reasonable way to classify the $CW$-spectra that realize this complex? 
 
\medskip
In this paper we take up these questions.  We also discuss how studying Floer theory from a homotopy perspective has been done in recent years, and how it has been applied with dramatic success. We state immediately that the applications that we  discuss in this paper are purely the choice of the author.  There are many other fascinating applications and advances that all help to make this an active and exciting area of research. We apologize in advance to researchers  whose work we will not have the time or space to discuss.

This paper is organized as follows.  In section 2 we discuss the three questions raised above.  We give a new take on how these questions were originally addressed in \cite{cjs}, and give some of the early applications of this perspective.    In section 3 we describe how Lipshitz and Sarkar used this perspective to define the notion of ``Khovanov homotopy" \cite{lipsark}\cite{lipsark2}.  This is a stable homotopy theoretic realization of Khovanov's homology theory \cite{khov} which in turn is a categorification of the Jones polynomial invariant of knots and links.  In particular we describe how the homotopy perspective Lipshitz and Sarkar used give more subtle and delicate invariants than the homology theory alone, and how these invariants have been applied.    In section 3 we describe Manolescu's work on an equivariant  stable homotopy theoretic view of Seiberg-Witten Floer theory.  In the case of his early work \cite{Man03}\cite{Man07}, the group acting is the circle group $S^1$.  In his more recent work \cite{Man13} \cite{Man16} he studied $Pin(2)$-equivariant Floer homotopy theory and used it to give a dramatic solution (in the negative) to the longstanding question about whether all topological manifolds admit triangulations.   In section 4  we describe the Floer homotopy theoretic methods of Kragh \cite{kragh} and of Abouzaid-Kragh \cite{Abou-Kra} in the study of the symplectic topology of the cotangent bundle of a closed manifold, and how they were useful in studying Lagrangian immersions and embeddings inside the cotangent bundle.

\paragraph{Acknowledgments.} The  author would like to thank M. Abouzaid, R. Lipshitz, and C. Manolescu for helpful comments on a previous draft of this paper.

\section{The homotopy theoretic foundations}
\subsection{Realizing chain complexes} 
We begin with a purely homotopy theoretic question:  Given a finite chain complex, $C_*$, 
$$
C_n \xr{\p_n} C_{n-1} \xr{\p_{n-1}} \cdots \xr{\p_2}C_1 \xr{\p_1}C_0
$$
how can one classify the finite $CW$-spectra $\bf X$ whose associated cellular chain complex is $C_*$?  This is question {\bf Q3} above.  

Of course one does not need to restrict  this question to finite complexes, but that is a good place to start.  In particular it is motivated
by Morse theory, where, given a Morse function on a closed, $n$-dimensional  Riemannian manifold, $f : M^n \to \br$, one has a  corresponding   ``Morse - Smale" chain complex $C_*^f$ 

$$
C^f_n \xr{\p_n} C^f_{n-1} \xr{\p_{n-1}} \cdots \xr{\p_2}C^f_1 \xr{\p_1}C^f_0
$$
Here $C^f_p$ is the free abelian group generated by the critical points of $f$ of index $p$.  The boundary homomorphisms are as described above  (\ref{morsecomp}). Of course in this case the Morse function together with the Riemannian metric define a $CW$ complex homotopy equivalent to $X$. In Floer theory that is not the case.   One does start with geometric information that allows for the definition of a chain complex, but knowing if this complex comes, in a natural way from the data of  the Floer theory, is not at all clear, and  was the central question of  study in \cite{cjs}.  

This homotopy theoretic question was addressed more specifically in \cite{floeroslo}.  Of central importance in this study was to understand how the attaching maps  of the cells in a finite $CW$-spectrum can be understood geometrically, via the theory of (framed) cobordism of manifolds with corners.  We will review these ideas in this section and recall some basic examples of how they can be applied. 

\med
 By assumption, the chain complex $C_*$ is finite, so  each $C_i$ is a finitely generated free abelian group.  Let   $\cb_i$ be a basis for $C_i$.       Let $\bs$ denote the sphere spectrum.  For each $i$,   consider the free $\bs$-module spectrum generated by $\cb_i$,   $$
\ce_i = \bigvee_{\alpha \in \cb_i}  \bs. $$    There is a natural isomorphism  
$$
  H_0(  \ce_i)  \simeq C_i.
$$

\begin{definition}\label{realize}  We say that a finite spectrum ${\bf X} $  realizes  the complex $C_* $,  if there exists a filtration
 of spectra converging to $X$,
 $$
 X_0 \hk X_1 \hk \cdots \hk X_n =  {\bf X}
 $$
satisfying the following properties:
\begin{enumerate}
\item  There is an equivalence of   the subquotients
$$
X_i/X_{i-1} \simeq \Sigma^i \ce_i,
$$ 
%and
\item
The induced composition  map in  integral homology,
\begin{align}
  &\tilde H_i(X_i/X_{i-1}) \xr{\delta_i} \tilde H_{i}(\Sigma X_{i-1}) \xr{\rho_{i-1}} H_{i}(\Sigma (X_{i-1}/X_{i-2}))   \notag \\
  &= H_0(\ce_i)      \to H_0(\ce_{i-1}) \notag \\
&= C_i  \to  C_{i-1} \notag
\end{align}
is the boundary homomorphism, $\p_i  $.     \notag
\end{enumerate}
Here the    ``subquotient" $X_i/X_{i-1}$ refers to the homotopy cofiber of the map $X_{i-1} \to X_i$, the map $\rho_i : X_i \to X_i/X_{i-1}$ is the projection map, and  the map $\delta_i : X_i /X_{i-1} \to \Sigma X_{i-1}$ is the Puppe extension of the homotopy cofibration sequence $X_{i-1} \to X_i \xr{\rho_i} X_i/X_{i-1}$.
\end{definition}

To classify finite spectra that realize  a given finite complex, in \cite{cjs} the authors    
 introduced  a category $\cj $ whose objects are the nonnegative  integers $\bz^+$,  and whose non-identity morphisms from $i$ to $j$  is empty for $i  \leq j$,   for $i > j+1$ it  is defined to be the    one point compactification,
 $$
 Mor_{\cj}(i, j) \cong (\br_+)^{i-j-1}\cup \infty
 $$
 and $Mor_{\cj}(j+1, j) $ is defined to be the two point space, $S^0$.
 Here  $\br_+$ is the space of nonnegative real numbers.   Composition in this category can be viewed in the following way.  Notice that  for $i > j+1$ $Mor_{\cj}(i,j)$ can be viewed as the one point compactification of the space $J(i,j)$ consisting of sequences of real numbers $\{\lambda_k\}_{k\in \bz}$ such that
 
 \begin{align}
 \lambda_k &\geq 0   \quad \text{for all} \, \, k \notag \\
 \lambda_k &= 0 \quad \text{unless}  \,  i > k > j. \notag
 \end{align}

 For consistency of notation we write $Mor_{\cj}(i,j) = J(i,j)^+$.  
 Composition of morphisms $ J(i,j)^+ \wedge J(j,k)^+ \to  J(i,k)^+ $ is then induced by addition of sequences. In this smash product the basepoint is taken to be $\infty$.  Notice that this map is basepoint preserving.   Given integers $p > q$, then there are subcategories  $\cj^p_q$ defined to be  the full subcategory generated by integers $q \geq m \geq p.$  The category $\cj_q$ is the full subcategory of $\cj$ generated by all integers $m \geq q$.

  The following   is a recasting of a discussion   in \cite{cjs}.

 \begin{theorem}\label{realmod}  The realizations of the  chain complex $C_*$ by    finite spectra   correspond to extensions of the association
  $j \to \ce_j$  to basepoint preserving functors $Z : \cj_0 \to Spectra$,  with the property that for each $j \geq 0$, the map obtained by the application of morphisms, 
\begin{align}
   Z_{j+1, j} : J(j+1,j)^+ &\wedge \ce_{j+1} \to \ce_{j}  \notag  \\
  S^0 &\wedge \ce_{j+1} \to \ce_{j} \notag
  \end{align}
   induces the boundary homomorphism $\p_{j+1} $ on the level of homology groups.   Here $Spectra$ is  a symmetric monoidal category of spectra (e.g the category of symmetric spectra), and by  a ``basepoint preserving functor" we mean one that maps $\infty$ in $J (p,q)^+$ to the constant map   $Z(p) \to Z(q)$.  
  \end{theorem}

  \med
  Since this is a recasting of a result in \cite{cjs} we supply a proof here.
  
  \begin{proof}  Suppose one has a functor $Z : \cj_0 \to Spectra$ satisfying the properties described in the theorem.  
  One defines a ``geometric realization" of this functor as follows.

          As described in \cite{cjs},  given   a functor to the category of spaces,
 $  Z : \cj_q \to  Spaces_*$, where $Spaces_*$ denotes the category of based topological spaces,  one can take its geometric realization,
\begin{equation}\label{zrealize}
 |Z| = \coprod_{q\leq j  } Z(j)\wedge J(j, q-1)^+ / \sim
\end{equation}
 where one identifies the image of $Z(j) \wedge  J(j,i)^+\wedge J(i, q-1)^+$  in $Z(j) \wedge J(j, q-1)^+$  given by composition of morphisms, with its image in $Z(i) \wedge J(i, q-1)$ defined by application of morphisms $Z(j) \wedge J(j,i)^+ \to Z(i)$. 

 For a functor whose value is in $Spectra$, we replace the above construction  by a coequalizer, in the following way:

 Let $Z : \cj_q \to Spectra$.  Define two maps of spectra,
 \begin{equation}\label{iotamu}
 \iota, \mu :  \bigvee_{q \leq j} Z(j) \wedge  J(j,i)^+\wedge J(i, q-1)^+  \la  \bigvee_{q\leq j  } Z(j)\wedge J(j, q-1)^+.
\end{equation}
 The first map $\iota$ is induced by the composition of morphisms in $\cj_q$,  $J(j,i)^+ \wedge  J(i, q-1)^+ \hk J(j, q-1)^+.$  The second map $\mu$ is the given by the wedge of maps,
 $$
  Z(j) \wedge J(j,i)^+\wedge J(i, q-1)^+  \xr{\mu_q \wedge 1} Z(i)\wedge J(i, q-1)^+
 $$ where $\mu_q : Z(j) \wedge J(j,i)^+ \to Z(i)$ is given by application of morphisms.

 \med
 \begin{definition}\label{georeal}  Given a functor $Z:  \cj_q \to Spectra$ we define its geometric realization $|Z|$ to be the  homotopy  coequalizer (in the category $Spectra$) of the two maps,
 $$
  \iota, \mu :  \bigvee_{q \leq j} \bigvee_{i = q}^{j-1}Z(j) \wedge  J(j,i)^+\wedge J(i, q-1)^+  \la  \bigvee_{q\leq j  } Z(j)\wedge J(j, q-1)^+.
  $$
  \end{definition}
  
  Given a functor $Z : \cj_0 \to Spectra$ as in the statement of the theorem, we claim the geometric realization $|Z|$ realizes the chain complex $C_*$.  In this realization, the filtration $$|Z|_0 \hk |Z|_1 \hk \cdots |Z|_k \hk |Z|_{k+1} \hk \cdots \hk |Z|_n = |Z| $$
   is the natural one where $|Z|_k $ is the homotopy coequalizer of $\iota$ and $\mu$ restricted to
\begin{equation}\label{zeekay}
   \iota, \, \mu :  \bigvee_{j=0}^k\bigvee_{i = 0}^{j-1} Z(j) \wedge  J(j,i)^+\wedge J(i, -1)^+  \la  \bigvee_{j=0}^k Z(j)\wedge J(j, -1)^+.
\end{equation}

To see that this is a filtration with the property that $|Z|_k/|Z|_{k-1} \simeq \Sigma^k \ce_k$, notice that    for a based space or spectrum, $Y$,  the smash product $J(n,m)^+\wedge Y$  is  homeomorphic  iterated cone,
\begin{equation}\label{cone}
  Y \wedge J(n,m)^+ \cong  c^{n-m-1}(Y).
\end{equation}
In particular $Z(k) \wedge J(k, -1)^+$ is the iterated cone $ c^k(Z(k)) = c^k(\ce_k).$  Now lets consider how $|Z|_k$ is built from $|Z|_{k-1}$.  Let $\p J(k, -1)$ denote the subspace of $J(k, -1)$ defined to consist of those sequences of real numbers $\{\lambda_i\}_{i\in \bz}$ such that
  \begin{align}
 \lambda_i &\geq 0   \quad \text{for all} \quad  i \notag \\
 \lambda_i &= 0 \quad \text{unless}  \quad  k > i > -1, \notag \\
 \text{there is at least one value of $i$ for}  \, &k > i > -1, \text{ such that} \, \lambda_i = 0.  \notag
 \end{align}

 It is easy to check that $\p J(k,-1)$ is  homeomorphic to $\br^{k-1}$, and its one-point compactification
 $\p J(k, -1)^+$ is homeomorphic to the sphere $S^{k-1}$.   
 
 Now consider the map $\hat \mu :   Z(k)  \wedge   \p J(k, -1)^+\to |Z|_{k-1}$ defined  by the maps
 $$\mu : Z(k) \wedge J(k, i)^+ \wedge J(i, -1)^+ \to |Z|_{k-1}
$$
as described in (\ref{iotamu}) above. In particular  $\mu = Z_{k,i} \wedge 1$ where $Z_{k,i} : Z(k) \wedge J(k, i)^+  \to    Z(i) $   is given by the application of morphisms.

It is then clear from the definition of $|Z|_k$  (\ref{zeekay})
  that it is the homotopy cofiber of the map

\begin{align}
 \hat \mu : \p J(k, -1)^+ \wedge Z(k) &\to |Z|_{k-1}  \notag \\
 \cong S^{k-1} \wedge Z(k) &\to |Z|_{k-1}.
\end{align}

Therefore the homotopy cofiber of the inclusion
$|Z|_{k-1} \hk |Z|_k$ is $S^k \wedge Z(k)$, which by hypothesis is $S^k \wedge \ce_k$.

 \med

    We now show  that the map
\begin{align}
   Z_{j+1, j} : J(j+1,j)^+ &\wedge \ce_{j+1} \to \ce_{j}  \notag  \\
  S^0 &\wedge \ce_{j+1} \to \ce_{j} \notag
  \end{align} induces  the boundary homomorphism $\p_j$ in homology.   To see this, recall that the boundary homomorphism
  is given by applying homology to the composition
  \begin{align}
 \p :  |Z|_{j+1}/|Z|_j &\xr{\alpha}  \Sigma |Z|_j  \xr{\rho}   \Sigma |Z|_j/|Z|_{j-1} \notag \\
\p :  S^{j+1}\wedge Z(j+1)     &\to S^{j+1} \wedge Z(j) \notag  \\
 \p :   S^{j+1}\wedge \ce_{j+1}     &\to S^{j+1} \wedge \ce_j \notag 
   \end{align}
Here $\alpha :  |Z|_{j+1}/|Z|_j \to  \Sigma |Z|_j $ is the connecting map in the Puppe extension of the cofibration sequence
$$
|Z|_j \hk |Z|_{j+1} \to |Z|_{j+1}/|Z|_j 
$$ and $\rho : \Sigma |Z|_j  \to  \Sigma |Z|_j/|Z|_{j-1}$ is projection onto the (homotopy) cofiber.   By the above discussion we have a homotopy cofibration sequence
$$
\cdots \to Z(j+1) \wedge \p J(j+1, -1)  \xr{\hat \mu} |Z|_j \to |Z|_{j+1} \xr{\rho}  |Z|_{j+1}/ |Z|_j \xr{\alpha} \Sigma |Z|_j \to \cdots
$$
So we have a homotopy commutative diagram

$$
\begin{CD}
|Z|_{j+1}/ |Z|_j @>\alpha >> \Sigma |Z|_j  \\
@V \simeq VV    @VV = V \\
\Sigma \p J(j+1, -1) \wedge Z(j+1)   @>> \Sigma \hat \mu >  \Sigma |Z|_j 
\end{CD}
$$

So the boundary map $\p$ given above can be viewed as the composition
\begin{equation} \label{comp} \p :  \Sigma \p J(j+1, -1) \wedge Z(j+1)  \xr{ \Sigma \hat \mu }   \Sigma |Z|_j  \xr{\rho} \Sigma |Z|_j/|Z|_{j-1}. \end{equation} 

Notice that when restricted to any $\Sigma J(j+1, r)^+ \wedge J(r, -1)^+ \wedge Z(j+1)$ for $r < j$, the composition $\p$  is naturally null homotopic.   This says that  $\p$  factors through the umkehr map (Pontrjagin-Thom construction)
\begin{equation} \label{pontthom}
 \p :  \Sigma \p J(j+1, -1) \wedge Z(j+1)  \xr{\tau} \Sigma J(j+1, j)^+ \wedge int(J(j, -1)^+ \wedge Z(j+1) \xr{ \bar \p} \Sigma |Z|_j/|Z|_{j-1}.
 \end{equation} 
 Here for $r > s$
$int (J(r, s)) $ is  the interior of the space $J(r,s)$.   This consists of sequences 
$\{\lambda_i\}_{i\in \bz}$ such that
  \begin{align}
 \lambda_i &\geq 0   \quad \text{for all} \quad  i \notag \\
 \lambda_i & = 0 \quad \text{unless}  \quad  r > i > s   \quad \text{in which case} \quad  \lambda_i > 0.  \notag 
  \end{align}  In particular $int (J(r, s)) \cong \br^{r-s-1}$ and so $int (J(r, s))^+ \cong S^{r-s-1}$.  
  
  Thus the boundary map $\p$ can be viewed as the composition
\begin{align}
 \p :  \Sigma \p J(j+1, -1) \wedge Z(j+1)  &\xr{\tau} \Sigma J(j+1, j)^+ \wedge int(J(j, -1)^+ \wedge Z(j+1) \xr{ \bar \p} \Sigma |Z|_j/|Z|_{j-1} \notag \\
 \p : S^{j+1}  \wedge Z(j+1) &\xr{\tau}   J(j+1,j) \wedge S^{j+1}   \xr{ \bar \p} \Sigma |Z|_j/|Z|_{j-1} \simeq S^{j+1}\wedge Z(j) \notag
 \end{align}
 where $\tau$ is an equivalence. But by the definition of $\bar \p$   (\ref{pontthom}) in terms of $ \mu = Z_{j, j-1} : J(j+1, j)^+ \wedge Z(j) \to Z(j-1)$,  we see that the boundary homomorphism $\p_j : C_j \to C_{j-1} $ is  defined by the application of morphisms, as asserted.  
  \end{proof}
  
So  a functor $Z : \cj_0 \to Spectra$ satisfying the properties specified in Theorem \ref{realmod}   defines a geometric realization $|Z|$ which is a finite spectrum.   Consider how the data of the functor $Z$ defines the $CW$-structure of $|Z|$.     Clearly $|Z|$ will have one cell of dimension $i$ for every element of $\pi_0(Z(i)) = \pi_0(\ce_i) = \cb_i$.  The attaching maps were described in \cite{cjs}, \cite{cotangent}, \cite{floeroslo} in the following way.

In general assume that  $X$ be a finite $CW$-spectrum with skeletal filtration
 $$
X_0 \hk X_1 \hk \cdots \hk X_n = X.
$$
In particular each map   $X_{i-1} \hk X_i$ is a cofibration, and we call its cofiber  $K_i = X_i /(X_{i-1})$.  This is a wedge of (suspension spectra) of spheres of dimension $i$.    
$$
K_i \simeq \bigvee_{\cd_i} \Sigma^i\bs
$$
where $\cd_i$ is a finite indexing set.

 As was discussed in \cite{cjs} one can then   ``rebuild" the homotopy type of the $n$-fold suspension, $\Sigma^n X,$
as the union of iterated cones and suspensions of the $K_i$'s,
\begin{equation}\label{decomp}
\Sigma^n X \simeq \Sigma^nK_0 \cup c(\Sigma^{n-1} K_1) \cup \cdots \cup c^i(\Sigma^{n-i} K_i) \cup \cdots \cup c^n K_n.
\end{equation}

This decomposition can be described as follows.
Define a map $\delta_i : \Sigma^{n-i}K_i \to \Sigma^{n-i+1}K_{i-1}$ to be the iterated suspension of the composition,
 $$
 \delta_i : K_i \to \Sigma X_{i-1} \to \Sigma K_{i-1}
 $$
 where the two maps in this composition come from the cofibration sequence, $X_{i-1}\to X_i \to K_i \to \Sigma X_{i-1} \cdots$.     As was pointed out in \cite{cotangent}, this induces a ``homotopy chain complex",

  \begin{equation}\label{htpychain}
K_n \xr{\delta_n} \Sigma K_{n-1} \xr{\delta_{n-1}} \cdots \xr{\delta_{i+1}}\Sigma^{n-i }K_{i } \xr{\delta_{i}} \Sigma^{n-i+1}K_{i-1} \xr{\delta_{i-1}} \cdots \xr{\delta_1}\Sigma^{n}K_0 = \Sigma^n X_0.
\end{equation}
 We refer to this as a homotopy chain complex  because examination of the defining cofibrations leads to canonical null homotopies of the compositions,
 $$
\delta_j \circ \delta_{j+1}.
$$
 This canonical null homotopy
defines an extension of $\delta_j$ to the    mapping cone of $\delta_{j+1}$:
$$
c(\Sigma^{n-j-1}K_{j+1}) \cup_{\delta_{j+1}} \Sigma^{n-j}K_j   \la \Sigma^{n-j+1}K_{j-1}.
$$
More generally,  for every $q$, using these null homotopies,   we have an extension  to the iterated mapping cone,
\begin{equation}\label{attach}
c^q(\Sigma^{n-j-q}K_{j+q}) \cup c^{q-1}(\Sigma^{n-j-q+1}K_{j+q-1}) \cup \cdots \cup c(\Sigma^{n-j-1}K_{j+1}) \cup_{\delta_{j+1}} \Sigma^{n-j}K_j  \la  \Sigma^{n-j+1}K_{j-1}.
\end{equation}

In other words, for each $p > q$, these null homotopies define  maps of spectra,
\begin{equation}\label{phi}
\phi_{p,q} :  c^{p- q-1}\Sigma^{n-p}K_p \to \Sigma^{n-q}K_q.
\end{equation}
The cell attaching data in the $CW$-spectrum $X$  as in (\ref{decomp}) is then defined via the maps $\phi_{p,q}$.  

Given a functor $Z : \cj_0 \to Spectra$ satisfying the hypotheses of Theorem \ref{realmod}, We have that 
$$
K_p = |Z|_p/|Z|_{p-1} \simeq \bigvee_{\cb_p} \Sigma^p \bs
$$
and the attaching maps 
\begin{align}
\phi_{p,q} : &c^{p- q-1}\Sigma^{-p}K_p \to \Sigma^{-q}K_q  \\
&J(p,q)^+ \wedge\Sigma^{-p} |Z|_p/|Z|_{p-1}     \to \Sigma^{-q}|Z|_q/|Z|_{q-1} \notag \\
&J(p,q)^+ \wedge Z(p) \to Z(q) \notag
\end{align}
is given by the application of the morphisms of the category $\cj_0$ to the value of the functor $Z$.

\subsection{Manifolds with corners and framed flow categories}

In order to make use of Theorem \ref{realmod}  in the setting of Morse and Floer theory, one needs a more geometric way of understanding the homotopy theoretic information contained in a functor $Z : \cj_0 \to Spectra$ satisfying the hypotheses of the theorem.   As is common in algebraic and differential topology, the translation between homotopy theoretic information and geometric information is done via cobordism theory and the Pontrjagin-Thom construction.  

In the case when the $CW$-structure comes from a Morse function $f : M^n \to \br$ on a closed $n$-dimensional Riemannian manifold,  the attaching maps $\phi_{p,q}$  defining a functor $Z_f : \cj_0 \to Spectra$,  were shown in \cite{cotangent}, \cite{floeroslo} to come from the cobordism-type of the moduli spaces of gradient flow lines connecting critical points of index $p$ to those of index $q$.  The relevant cobordism theory is \sl framed cobordism   of manifolds with corners.  \rm Thus to extend this idea to the Floer setting,  with the goal of developing a ``Floer homotopy type", one needs to understand certain cobordism-theoretic properties of the corresponding moduli spaces of gradient flows of the particular Floer theory. (In \cite{floeroslo} the author also considered how other cobordism theories give rise not to ``Floer homotopy types" but rather to ``Floer module spectra", or said another way, ``Floer $E$-theory" where $E_\bullet $  is a generalized homology theory.    We refer the reader to \cite{floeroslo} for details.)   In order to make this idea precise we recall some basic facts about cobordisms of manifolds with corners.  The main reference for this is Laures's paper \cite{laures}.

 Recall that an $n$-dimensional manifold with corners, $M$,  has charts which are local homeomorphisms with  $\br_+^n$.  Here $\br_+$ denotes the nonnegative real numbers and $\br_+^n$ is the $n$-fold cartesian product of $\br_+$.    Let $\psi : U \to  \br_+^n$ be a chart of a manifold with corners $M$.  For $x \in U$,  the number of zeros of this chart, $c(x)$ is independent of the chart.   One defines a \sl face \rm of $M$ to be a connected component of the space  $\{m{\in}M \, \text{such that} \, c(m)=1\}$.

 Given an integer $k$, there is a notion of a manifold with corners having ``codimension k",  or a $\langle k\rangle$-manifold.  We recall the definition from \cite{laures}.

 \begin{definition}\label{kmanifold} A $\langle k\rangle$-manifold  is a manifold  with corners,   $M$,  together with an ordered $k$-tuple  $(\partial{_1}M,...,\partial{_k}M)$ of unions of faces of $M$  satisfying the following properties.
 \begin{enumerate}
 \item Each $m{\in}M$ belongs to $c(m)$ faces
 \item $\partial{_1}M  \, {\cup}  {\cdots}  {\cup} \,  \partial_{k}  M = \partial {M} $
 \item For all $1{\leq}i{\neq}j{\leq}k$, $\partial{_i}M \, {\cap}  \, \partial{_j}M$ is a face of both  $\partial{_i}M$ and $\partial{_j}M$.
 \end{enumerate}

\end{definition}

  The  archetypical example of a  $\langle k \rangle$-manifold is $ \br_+^k$.  In this case
 the face $F_j \subset \br_+^k$ consists of those $k$-tuples with the $j^{th}$-coordinate equal to zero.

As described in \cite{laures}, the data of a $\langle k \rangle$-manifold can be encoded in a   categorical way  as follows.  Let $\underbar{2}$ be the partially ordered set with two objects, $\{0, 1\}$, generated by a single nonidentity morphism $0 \to 1$.  Let $\ut^k$ be the product of $k$-copies of the category $\ut$.  A $\langle k \rangle$-manfold $M$ then defines a functor from $\ut^k$   to the category of topological spaces,  where for an object  $a = (a_1, \cdots , a_k) \in \ut^k$,
$M(a)$ is the intersection of the faces $\p_i(M)$   with $a_i = 0$.  Such a functor is a $k$-dimensional cubical diagram of spaces, which, following Laures' terminology, we refer to as a $\langle k \rangle$-diagram, or a $\langle k \rangle$-space.    Notice that $\br_+^k(a) \subset \br_+^k$ consists of those $k$-tuples of nonnegative real numbers so that the $i^{th}$-coordinate is zero for every $i$  with $a_i=0$.  More generally, consider the $\langle k\rangle$-Euclidean space, $\br_+^k \times \br^n$,  where the value on $a \in \ut^k$ is $\br_+^k(a) \times \br^n$.      In general  we refer to a functor $\phi : \ut^k \to \cc$ as a $\langle k\rangle$-object in the category $\cc$.

    In this section we will consider embeddings of manifolds with corners into Euclidean spaces $M \hk \br_+^k \times \br^n$ of the form given by the following definition.

\begin{definition}\label{embed} A ``neat embedding" of a $\langle k \rangle$-manifold $M$ into   $ \br_+^k \times \br^m  $  is a natural transformation of $\langle k \rangle$-diagrams
$$
e : M \hk  \br_+^k \times \br^m
$$   that satisfies the following properties:
\begin{enumerate}
\item  For each $a{\in}\underline{2}^k$, $e(a)$ is an embedding.
\item
For all $b < a$, the intersection   $M(a) \cap  \left( \mathbb{R}^k_+(b)  \times  \mathbb{R}^m\right) = M(b)$,  and this intersection is perpendicular.  That is,    there is some $\epsilon > 0$ such that $$M(a) \cap  \left(\mathbb{R}^k_+(b)  \times  [0,\epsilon)^k(a-b) \times \mathbb{R}^m\right) = M(b)  \times [0,\epsilon)^k(a-b).$$
\end{enumerate}
Here $a-b$ denotes the object of $\ut^k$ obtained by subtracting the $k$-vector $b$ from the $k$-vector $a$.
\end{definition}

\med
In \cite{laures} it was proved   that every $\langle k \rangle$-manifold neatly embeds in $\mathbb{R}^{k}_+{\times} \br^N$ for $N$ sufficiently large.  In fact it was proved there that a manifold with corners, $M$,  admits  a neat embedding into $\mathbb{R}^{k}_+{\times} \br^N$ \sl if and only if \rm $M$ has the structure of a $\langle k\rangle$-manifold.  Furthermore in \cite{genauer} it  is shown that the connectivity of the space of neat embeddings, $Emb_{\langle k\rangle}(M; \mathbb{R}^{k}_+{\times} \br^N)$ increases with the dimension $N$.

\med
Notice that    an embedding of manifolds with corners,  $e : M \hk   \br_+^k  \times  \br^m  $,    has a well defined normal bundle.  In particular, for any pair of objects in $\ut^k$,  $a > b$,   the normal bundle of $e(a) : M(a) \hk  \br_+^k(a)  \times \br^m$, when restricted to $M(b)$, is the normal bundle of $e(b) : M(b) \hk  \br_+^k(b) \times \br^m$.   

\med
These embedding properties of $\langle k \rangle$ - manifolds make it clear that these are the appropriate manifolds to study for cobordism - theoretic information.   In particular, given an embedding $e: M \hk  \br_+^k \times \br^m$  the Thom space of the normal bundle, $Th (M, e)$,  has the structure of an $\langle k\rangle$-space, where for $a \in \ut^k$, $Th (M, e)(a)$ is the Thom space of the normal bundle of the associated embedding, $M(a) \hk \br_+^k (a) \times \br^N$.  We can then desuspend  and define the Thom spectrum, $M^\nu_e = \Sigma^{-N}Th (M,e)$,  to be the associated $\langle k\rangle$-spectrum. The Pontrjagin-Thom construction defines a map of $\langle k\rangle$-spaces, $$\tau_e : \left( \br_+^k \times \br^N\right)\cup \infty = (( \br_+^k)\cup \infty) \wedge S^N \to  Th (M, e).$$ Desuspending we get a map of $\langle k\rangle$-spectra, $\Sigma^\infty (( \br_+^k)\cup \infty) \to  M^\nu_e. $  Notice that  the homotopy type (as  $\langle k\rangle$-spectra) of $M^\nu_e$ is independent of the embedding $e$. We denote the homotopy type of this normal Thom spectrum as $M^\nu$, and  the Pontrjagin-Thom map,  $\tau : \Sigma^\infty (( \br_+^k)\cup \infty) \to M^\nu$.

\med 
Compact manifolds with corners, and in particular $\langle k \rangle$ - manifolds naturally occur as  the moduli spaces of flow lines of a Morse function, and in some cases, of a Floer function.  We first recall how they appear in Morse theory.

\med
Consider  a smooth, closed $n$-manifold $M^n$, and a smooth  Morse  function $f : M^n \to \br$. Given a Riemannian metric on $M$,  one studies the flow of the gradient vector field $\nabla f$.
 In particular a flow line is a curve $\gamma : \br \to M$ satisfying the ordinary differential
 equation,
 $$
 \frac{d}{dt}\gamma(s) + \nabla f (\gamma (s)) = 0.
 $$
 By the existence and uniqueness theorem for solutions to ODE's, one knows that if
 $x \in M$ is any point then  there is a unique flow line $\gamma_x$ satisfying $\gamma_x(0) = x$.  One then studies unstable and stable manifolds of the critical points,
 \begin{align}
 W^u(a) &= \{ x \in M \, : \, \lim_{t\to -\infty}\gamma_x (t) = a \}  \notag \\
 W^s(a) &= \{ x \in M \, : \, \lim_{t\to +\infty}\gamma_x (t) = a \}.   \notag
 \end{align}

The unstable manifold $W^u(a)$ is diffeomorphic to a disk $D^{\mu (a)}$, where $\mu (a)$ is the index of the critical point $a$.  Similarly    the stable manifold $W^s(a)$  is  diffeomorphic to a disk $D^{n-\mu (a)}$.

For a generic choice of Riemannian metric, the unstable manifolds and stable manifolds
intersect transversally, and their intersections,
$$
W(a,b) = W^u(a) \cap W^s (b)
$$
are smooth manifolds   of  dimension equal to the relative index,  $\mu (a) - \mu (b)$.  When the choice of metric satisfies these transversality properties,  the metric is  said to be ``Morse-Smale".     The manifolds $W(a,b)$ have   free $\br$-actions  defined by  ``going with the flow".  That is,   for $t \in \br$, and $x \in M$,
$$
t\cdot x = \gamma_x (t).
$$
The ``moduli space of flow lines" is the manifold
$$
\cm (a,b) = W(a,b) /\br
$$
and has dimension $\mu (a) - \mu (b) -1$.  These moduli spaces are not generally compact, but they have canonical compactifications which we now describe.

 In the case of  a Morse-Smale metric, (which we assume throughout the rest of this section),     there is a partial order on the finite set of critical points   given by    $a \geq b$ if $\cm(a,b) \neq \emptyset$.   We then define
\begin{equation} \label{compact}
\bar \cm (a,b) = \bigcup_{a=a_1 >a_2> \cdots > a_k = b} \cm(a_1, a_2) \times \cdots \times \cm (a_{k-1}, a_k),
\end{equation}

The topology of $\bar \cm (a,b)$ can be  described naturally, and is done so in many references including \cite{cjs1}.   $\bcm (a,b)$ is the space of   ``piecewise flow lines" emanating from $a$ and ending at $b$.  

 The following definition of a Morse function's  ``flow category " was also given in \cite{cjs1}.
  
\med
\begin{definition}  The \sl flow category \rm $\cc_f$  is a topological category associated to a Morse function $f : M \to \br$ where $M$ is a closed Riemannian manifold.  Its objects are  the critical points of $f$.    If $a$ and $b$ are two such critical points, then $Mor_{\cc_f}(a, b) = \bcm (a, b)$.   Composition is determined by the maps
$$
\bcm (a,b) \times \bcm (b, c) \hk \bcm (a, c)
$$
which are defined to be the natural embeddings into the boundary.
\end{definition}

The moduli spaces $\cm (a,b)$  have natural framings on their stable normal bundles (or equivalently, their stable tangent bundles)  that play an important role in this theory.  These framings are defined in the following manner.  Let $a > b$ be critical points.  Let $\eps > 0$ be chosen so that there are no critical values in the half open interval $[f(a)-\eps, f(a))$.  Define the \sl unstable sphere \rm to be the level set of the unstable manifold,
$$
 S^u(a) = W^u(a) \cap f^{-1}(f(a)-\eps).
 $$
The sphere $S^u(a)$  has dimension $ \mu(a) - 1$. Notice there is a natural diffeomorphism,
$$
\cm(a,b) \cong S^u(a) \cap W^s(b).
$$
This leads to the following diagram,
\begin{equation}\label{intersect}
\begin{CD}
W^s(b)   @>\hk >>  M \\
@A \cup AA   @AA\cup A \\
\cm (a,b)  @>>\hk >  S^u(a).
\end{CD}
\end{equation}
From this diagram one sees that the normal bundle $\nu$ of the embedding $\cm(a,b) \hk S^u(a)$ is the restriction of the normal bundle of $W^s(b) \hk M$.  Since $W^s(b)$ is a disk, and therefore contractible, this bundle is trivial. Indeed an orientation of $W^s(b)$  determines a homotopy class of trivialization, or a framing.  In fact  this framing determines a diffeomorphism of the bundle  to the product, $W^s(b) \times W^u(b)$.  Thus these orientations give the moduli spaces
$\cm(a,b)$   canonical normal framings, $\nu \cong \cm(a,b) \times W^u(b) . $

As was pointed out in \cite{cjs1}, these framings extend to the boundary of the compactifications, $\bcm (a, b)$.  In order to describe what it means for these framings to be ``coherent" in an appropriate sense, the following categorical approach was used in \cite{floeroslo}.  The  first step is  to abstract the basic properties of a flow category of a Morse functon.

\begin{definition}\label{cptcat} A smooth, compact  category  is a topological category $\cc$ whose objects form a discrete set, and whose whose morphism spaces, $Mor (a,b)$ are compact, smooth  $\langle k \rangle$ -manifolds, where $k = dim \, Mor(a,b)$.  The  composition maps,  $\nu : Mor (a,b) \times Mor (b, c) \to Mor(a, c)$,  are smooth codimension one embeddings (of manifolds with corners) whose images lie in  the boundary.  Moreover every point in the boundary of $Mor(a, c)$ is in the image under $\nu$ of a unique maximal sequence in $Mor(a, b_1) \times Mor (b_1, b_2) \times \cdots \times Mor(b_{k-1}, b_k)\times Mor (b_k, c)$ for some objects $\{b_1, \cdots , b_k\}$.

\med
A smooth, compact category  $\cc$  is said to be a ``Morse-Smale" category if the following
additional properties are satisfied.
\begin{enumerate}
\item The objects of $\cc$ are partially ordered by the  condition
$$
 a \geq  b   \quad \text{if} \quad Mor (a,b) \neq \emptyset.
$$
\item $Mor(a,a) = \{identity \}$.
\item  There is a set map, $\mu : Ob (\cc) \to \bz$, which preserves the partial ordering,
such that  if $a > b$,
$$
dim \, Mor(a,b) = \mu (a) - \mu (b) -1.
$$
The map $\mu$ is known as an ``index" map.   A Morse-Smale category such as this  is said to have finite type, if there are only finitely many objects of any given index, and  for each pair of objects  $a > b$,  there are only finitely many objects
$c$ with $a > c > b$.   For ease of notation we write $k(a,b) = \mu(a) -\mu (b) -1$.  
\end{enumerate}
\end{definition}

The following is a folk theorem that goes back to the work of Smale and Franks \cite{franks}  although a proof of this fact did not appear in the literature
until much later \cite{qin}.

\begin{proposition}    Let $f : M \to \br$ be  smooth Morse function on a closed Riemannian manifold with a Morse-Smale metric. Then  the  compactified moduli space of piecewise flow-lines, $\bar \cm (a,b)$ is a  smooth  $\langle  k(a,b) \rangle $ -  manifold.   
\end{proposition}  

Using this result,  as well as an associativity result for the gluing maps $\bcm(a,b) \times \bcm (b, c) \to  
\bcm (a,c)$ which was eventually proved in \cite{qin},  it was   proven in \cite{cjs}  that the flow category $\cc_f$ of such a Morse-Smale function is indeed a Morse-Smale smooth, compact category according to Definition \ref{cptcat}.   

\med
\noindent
\bf Remark.  \rm The fact that \cite{cjs1} was never submitted for publication was due to the fact that the ``folk theorem" mentioned above, as well as the associativity of gluing, both of which the authors of \cite{cjs1} assumed were ``well known to the experts", were indeed not in the literature, and their proofs which was eventually provided in \cite{qin}, were analytically more complicated than the authors imagined.

\med
 In order to define the notion of  ``coherent framings" of the moduli spaces $\bcm (a,b)$,  so that we may apply the Pontrjagin-Thom construction coherently,  we need to study an associated category, enriched in spectra, defined using the stable normal bundles of the moduli spaces of flows.

 \med
 \begin{definition}\label{normalthom} Let $\cc$ be  a smooth, compact category of finite type satisfying the Morse-Smale condition.   Then a ``normal Thom spectrum" of the category $\cc$  is a category, $\cc^\nu$, enriched over spectra, that satisfies the following properties.
 \begin{enumerate}
 \item The objects of $\cc^\nu$ are the same as the objects of $\cc$.
 \item The morphism spectra  $Mor_{\cc^\nu}(a,b)$ are $\langle k(a,b)\rangle$-spectra,  having  the homotopy type of the normal Thom spectra $Mor_{\cc}(a,b)^\nu$, as $\langle k(a,b)\rangle$-spectra. The composition maps, $$\circ : Mor_{\cc^\nu}(a,b)\wedge Mor_{\cc^\nu}(b,c) \to Mor_{\cc^\nu}(a,c)$$
 have the homotopy type of the maps, $$
 Mor_{\cc}(a,b)^\nu \wedge Mor_{\cc}(b,c)^\nu \to Mor_{\cc}(a,c)^\nu$$ of the  Thom spectra of the stable normal bundles corresponding to the composition maps in $\cc$, $Mor_{\cc} (a,b) \times Mor_{\cc} (b,c) \to Mor_{\cc}(a,c)$.  Recall that these   are maps of $\langle k(a,c)\rangle$-spaces induced by the inclusion of a component of the boundary.
 \item The morphism spectra are equipped with  ``Pontrjagin-Thom  maps"  $\tau_{a,b} : \Sigma^\infty (J(\mu (a), \mu (b))^+) = \Sigma^\infty ( (  \br_+^{k(a,b)})\cup \infty)) \to Mor_{\cc^\nu}(a,b)$    such that the following diagram  commutes:
 $$
 \begin{CD}
   \Sigma^\infty (J(\mu (a), \mu (b))^+) \wedge  \Sigma^\infty (J(\mu (b), \mu (c))^+)  @>>>  \Sigma^\infty (J(\mu (a), \mu (c))^+)  \\
   @V\tau_{a.b}\wedge \tau_{b,c} VV    @VV\tau_{a,c} V \\
   Mor_{\cc^\nu}(a,b) \wedge Mor_{\cc^\nu}(b,c)  @>>> Mor_{\cc^\nu}(a,c).
   \end{CD}
   $$ Here the top horizontal map is defined via the composition maps  in the category $\cj$,  and the bottom horizontal map is defined via the composition maps in $\cc^\nu$.
   \end{enumerate}
   \end{definition}
   
   With the notion of a ``normal Thom spectrum" of a flow category $\cc$,   the notion of a coherent $E^*$-orientation was defined in \cite{floeroslo}.  Here $E^*$ is a generalized cohomology theory represented by a commutative ($E_\infty$) ring spectrum $E$.  We recall that definition now.

\med   
First observe that  a  commutative ring spectrum $E$   induces a $\langle k\rangle$-diagram in the category of spectra (``$\langle k\rangle$-spectrum"),  $E\langle k\rangle$, defined in the following   manner.

 For $k = 1$, we let $E\langle1\rangle: \ut \to Spectra$  be defined by $E\langle 1\rangle(0) = S^0$, the sphere spectrum, and $E\langle1\rangle(1) = E$.  The image of the morphism $0 \to 1$ is the unit of the ring spectrum $S^0 \to E$.

 To define $E\langle k\rangle$ for general $k$,  let $a$ be an object of $\ut^k$.  We view $a$ as a vector of length $k$, whose coordinates are either zero or one.    Define $E\langle k\rangle (a)$ to be the multiple smash product of spectra, with a copy of $S^0$ for every every zero coordinate, and a copy of $E$ for every string of successive ones.  For example, if $k = 6$, and $a = (1, 0, 1,1, 0,1)$,  then $E\langle k\rangle (a) = E\wedge S^0\wedge E\wedge S^0\wedge E$.

 Given a morphism $a \to a'$ in $\ut^k$, one has a map $E\langle k\rangle (a) \to E\langle k\rangle (a')$ defined by combining the unit $S^0 \to E$ with the ring multiplication $E\wedge E \to E$.

Said another way, the functor $E\langle k\rangle : \ut^k \to Spectra$ is  defined by taking  the $k$-fold product functor $E\langle 1\rangle : \ut \to Spectra$ which sends $(0 \to 1)$ to $S^0 \to E$,  and then using the ring multiplication in $E$ to ``collapse" successive strings of $E$'s.

 \med
 This structure allows us to define one more construction.  Suppose $\cc$ is a smooth,  compact,  Morse-Smale category of finite type as in Definition \ref{cptcat}.  We can then define an associated category,  $E_{\cc}$,  whose objects are the same as the objects of $\cc$  and whose morphisms are given by the spectra,
 $$
 Mor_{E_{\cc}}(a,b)  = E\langle k(a,b)\rangle
 $$
 where $k(a,b) =  \mu(a) - \mu (b) -1$.  Here  $\mu (a)$ is the index of the object $a$ as in Definition \ref{cptcat}.  The composition law is the pairing,
\begin{align}
 E\langle k(a,b)\rangle  \wedge E\langle k(b,c)\rangle &=  E\langle k(a,b)\rangle \wedge S^0 \wedge  E\langle k(b,c)\rangle  \notag\\
 &\xr{1 \wedge u \wedge 1} E\langle k(a,b)\rangle \wedge E\langle 1\rangle \wedge  E\langle k(b,c)\rangle  \notag \\
 &\xr{\mu} E\langle k(a,c)\rangle. \notag
 \end{align}
  Here $u : S^0 \to E = E\langle 1\rangle$ is the unit.    This category encodes the multiplication in the ring spectrum $E$.

     \med
   \begin{definition}\label{eorient}
   An $E^*$-orientation of a   smooth, compact category of finite type satisfying the Morse-Smale condition,   $\cc$, is a functor, $u : \cc^\nu \to  E_{\cc}$, where $\cc^\nu$ is a normal Thom spectrum of $\cc$, such that on morphism spaces, the  induced map
   $$
   Mor_{\cc^\nu}(a,b) \to   E\langle k(a,b)\rangle
   $$
   is a map of $\langle k(a,b)\rangle$-spectra that defines an $E^*$ orientation of $ Mor_{\cc^\nu}(a,b) \simeq \bcm(a,b)^\nu$.
   \end{definition}

   The functor $u : \cc^\nu \to E_{\cc}$ should be thought of as a coherent family of  $E^*$- Thom classes for the normal bundles of the morphism spaces of $\cc$.   When $E = \bs$, the sphere spectrum, then an $E^*$-orientation, as defined here, defines a coherent family of framings of the morphism spaces, and  is equivalent to the notion of a \sl framing \rm of the category $\cc$, as defined in \cite{cjs}.

In \cite{cjs} the following was proved modulo the results of \cite{qin} which appeared much later.

\begin{theorem}\label{framed}  Let $f : M \to \br$ be a Morse function on a closed Riemannian manifold satisfying the Morse-Smale condition.   Then the flow category $\cc_f$ has a canonical structure as a  ``$\bs$-oriented, smooth, compact Morse-Smale category of finite type". That is, it is a  \sl ``framed, smooth compact Morse-Smale category".   The induced framings of the morphism manifolds $\bcm (a,b)$ are canonical extensions of the framings of the open moduli spaces $\cm (a,b)$ descrbed above (\ref{intersect}).
\end{theorem}

\med
The main use of the notion of compact, smooth, framed flow categories, is the following result.

\med
\begin{theorem}(\cite{cjs}, \cite{floeroslo})   Let $\cc$ be a compact, smooth, framed category of finite type satisfying the Morse-Smale property.  Then there is an associated, naturally defined functor $Z_{\cc} : \cj_0 \to Spectra$  whose geometric realization $|Z_{\cc}|$ realizes the associated ``Floer complex"
$$
\to \cdots \to C_{i+1} \xr{ \p_i }  C_{i} \xr{\p_{i-1}} \cdots \xr{\p_0} C_0.
$$
Here $C_j$ is the free abelian group generated by the objects $ a \in Ob (\cc)$ with $\mu (a) = j$, and the boundary homomorphisms are defined by
$$
 \p_{j-1} ([a])  =  \sum_{\mu (b) = j-1}  \#Mor(a,b) \, \cdot \, [b]
$$
where $\#Mor(a,b)$ is the framed cobordism type of the compact  framed manifold $Mor (a,b)$ which zero dimensional (since its dimension $= \mu (a) - \mu (b) -1 = 0$).  Therefore this cobordism type is simply an integer, which can be viewed as the oriented count of the number of points in $Mor(a,b)$.  
\end{theorem}

\med
\begin{proof} \sl Sketch. \rm 
The proof of this is sketched in \cite{cjs} and is carried out  in \cite{floeroslo} in the setting of $E^*$-oriented compact, smooth categories, for $E$ any ring spectrum.  The idea for defining the functor $Z_\cc$ is to use the Pontrjagin-Thom construction in the setting of framed manifolds with corners (more specifically, framed $\langle k \rangle$-manifolds).  Namely, one defines
$$
Z_{\cc}(j) = \bigvee_{a \in Obj (\cc) \, : \, \mu (a) = j} \bs.
$$
On the level of morphisms one needs to define, for every $i > j$,     a map of spectra
$$
Z_{\cc} (i,j) : J(i,j)^+ \wedge Z_{\cc} (i) \to Z_{\cc}(j).
$$
This is defined to be the wedge, taken over all $a \in Obj (\cc)$ with $\mu (a) = i$, and $b \in Obj (\cc)$ with $\mu (b) = j$,
of the maps 
$$
Z_{\cc}(a,b) : J(\mu(a), \mu (b))^+\wedge \bs_a \to \bs_b
$$
defined to be the  composition
\begin{equation}
  \Sigma^\infty (J(\mu(a), \mu (b)^+) \xr{\tau_{a,b} } Mor_{\cc^\nu} (a,b) \xr{\ u} \bs
  \end{equation}
%\end{equation}
%adjoint of the composition
%\begin{equation}
  %\Sigma^\infty (J(\mu(a), \mu (b)^+) \xr{\tau_{a,b} } Mor_{\cc^\nu} (a,b) \xr{\tilde u} \Sigma^\infty Mor_{\cc}(a, b)
%\end{equation}
where $\bs_a$ and $\bs_b$ are copies of the sphere spectrum indexed by $a$ and $b$ respectively in the definition of $Z_{\cc}(i)$ and $Z_{\cc}(j)$. $\tau_{a,b}$ is the Pontrjagin-Thom map, and  $ u$ is  $\bs$-normal orientation class (framing).  %the  equivalence of the Thom spectrum and the suspension spectrum defined by the framing Thom class $u : Mor_{\cc^\nu} (a,b) \to \bs$. 
Details can be found in \cite{floeroslo}.  \end{proof}

\med
Thus to define a \sl ``Floer homotopy type" \rm one is looking for a compact, smooth, framed category of finite type satisfying the Morse-Smale property.    The compact framed manifolds with corners   that constitute the morphism spaces, define, via the Pontrjagin-Thom construction, the attaching maps of the $CW$-spectrum defining this (stable) homotopy type.  

 We see from Theorem \ref{framed} that given a Morse-Smale function $f : M \to \br$, then its flow category satisfies these properties, and  it was proved in \cite{cjs}  that, not surprisingly, its Floer homotopy type is the suspension spectrum $\Sigma^\infty (M_+)$. The $CW$-structure is the classical one coming from Morse theory, with one cell of dimension $k$ for each critical point of index $k$. The fact that the compactified moduli spaces of flow lines,  which constitute the morphism spaces in this category,     together with their structure as framed manifolds with corners,  define the attaching maps in the $CW$-structure of $\Sigma^\infty M$,  can be viewed as a generalization of the well-known work of Franks in \cite{franks}.    

\med
As pointed out in \cite{cjs}, a distinguishing feature in the flow category of a Morse-Smale function $f : M \to \br$ is that the framing is \sl canonical. \rm See (\ref{intersect}) above. As also was pointed out in  \cite{cjs}, if one chooses a \sl different \rm framing of the flow category $\cc_f$,  then the ``difference" between the new framing and the canonical framing defines a functor
\begin{equation}\label{frame}
\Phi : \cc_f \to \cg L_1(\bs)
\end{equation}
where $ \cg L_1(\bs)$ is the category corresponding to the group-like monoid $GL_1(\bs)$ of ``units" of the sphere spectrum (see \cite{units}).  This monoid has the homotopy time of the colimit
$$
GL_1(\bs) \simeq \colim_{n \to \infty} \Omega^n_{\pm 1}S^n.
$$
Here the subscript denotes the path components of $\Omega^nS^n$ consisting of based self maps of the sphere $S^n$ of degree $\pm 1$.  By a minor abuse of notations, we let  $\Phi$ denote the framing of $\cc_f$  that defines the map (\ref{frame}). 

Passing to the geometric realizations of these categories, one gets a map
$$
\phi : M \to BGL_1(\bs).
$$
which we think of as the isomorphism type of a spherical fibration over $M$.  The following is also a result of \cite{cjs}.

\med
\begin{proposition}  If  $f : M \to \br$ is a Morse-Smale function on a closed Riemannian manifold, and its flow category $\cc_f$ is given a framing $\Phi$, then the \sl Floer homotopy type \rm of $(\cc_f , \Phi)$ viewed as a compact, smooth, framed category, is the Thom spectrum  of the corresponding stable spherical fibration, $M^\phi$ 
\end{proposition}

\med

\subsection{ The free loop space of a symplectic manifold and Symplectic Floer theory }

Let $(N^{2n}, \omega)$ be a symplectic manifold.  Here  $\omega$ is a closed, nondegenerate, skew symmetric bilinear form on the tangent bundle, $TN$.  Let $LN$ be its free loop space.  One of the earliest applications of Floer theory \cite{floerFix}  was to the (perturbed) symplectic action functional on $LN$.

Let $L_0N \subset LN$ be the path component of contractible loops in $N$.  Let $\tilde {L_0N }\xr{p} L_0N$ be its universal cover.  Explicitly,
$$
\tilde {LN} = \{(\gamma, \theta) \in L_0N \times Map(D^2, N) \, : \, \text{the restriction of $\theta$ to $S^1 = \p D^2$ is equal to} 
\quad \gamma  \} / \sim
$$
where the equivalence relation is given by $(\gamma, \theta_1) \sim (\gamma, \theta_2)$ if, when we combine $\theta_1$ and $\theta_2$ to define a map of the $2$-sphere,
$$
\theta_{1,2} = \theta_1 \cup \theta_2 : D^2 \cup_{S^1} D^2 = S^2 \to  N
$$
then $\theta_{1,2}$ is null homotopic.  In other words, $(\gamma, \theta_1)$  is equivalent to $ (\gamma, \theta_2)$ if $\theta_1$ and $\theta_2$ are homotopic maps relative to the boundary.  

One can then define the ``symplectic action" functional, 
\begin{align}\label{sympact}
\ca &: \tilde{L_0N} \to \br  \\
(\gamma, \theta) &\to \int_{D^2} \theta^*(\omega) \notag \\
\end{align}

The symplectic action descends to define an $\br/\bz$-valued function
$$
\ca : L_0N \to \br/\bz.
$$
One needs to perturb this functional by use of a Hamiltonian vector field in order to achieve nondegeneracy of critical points.  A Morse-type complex, generated by critical points, is then studied, and the resulting symplectic Floer homology   
has proved to be an important invariant.

We describe the situation when $N = T^*M^n$, the cotangent bundle of a closed, $n$-dimensional manifold, in more detail.  In particular this is a situation where one has a corresponding ``Floer homotopy type", defined via a compact, smooth flow category as described above.  This was studied by the author in \cite{cotangent}, making heavy use of the analysis of Abbondandolo and Schwarz in \cite{AS}.

Coming from classical mechanics, the cotangent bundle of a smooth manifold has a canonical symplectic structure.  It is defined as follows. 

Let $p : T^*M \to M$ be the projection map.  For $x \in M$ and $ u\in T^*_xM$,  let $z = (y, u) \in T^*M$, and consider the composition
 $$
 \theta _z: T_{z}(T^*M) \xr{dp}T_yM \xr{u}\br.
 $$
 
   This defines a $1$-form $\theta$ on $T^*M$, called the ``Liouville" $1$-form, and the symplectic form $\omega$ is defined to be the exterior derivative $\omega = d\theta$.  It is easy to check that $\omega$ is nondegenerate.    

 Let 
 $$
 H : \br/\bz \times T^*M \to \br
 $$
be a smooth function.  Such a map is called a ``time-dependent periodic Hamiltonian".  Using the nondegeneracy of the symplectic form, this allows one to define the corresponding ``Hamiltonian vector field" $X_H$ by requiring it to satisfy the equation 
 $$
 \omega (X_H(t,z), v) = -dH_{t,z} (v)
 $$
 for all $t \in \br/\bz$, $ z \in T^*M$, and $v \in T_{z}(T^*M)$. 
 We will be considering the space of $1$-periodic solutions, $\cp (H)$, of the Hamiltonian equation
 $$
 \frac{dz}{dt} =  X_H(t, z(t))
 $$ where $z : \br/\bz \to T^*M$ is a smooth function.

  \med
  Using a periodic time-dependent Hamiltonian one can define the perturbed 
  symplectic action functional 
 \begin{align}
 \ca_H : L (T^*M)   &\to \br \notag \\
z&\to \int z^*(\theta -Hdt) = \int_0^1(\theta (\frac{dz}{dt}) - H(t, z(t)) dt.
\end{align}
This is a smooth functional, and its   critical points  are the  periodic orbits of the Hamiltonian vector field,  $\cp (H)$.  
Now let $J$ be a $1$-periodic, smooth almost complex structure on $T^*M$, so that for each $t \in \br/\bz$, 
$$
\langle \zeta, \xi \rangle_{J_t} = \omega (\zeta, J(t,z)\xi ), \quad \zeta, \xi \in T_{z}T^*M, \, z \in T^*M,
$$
is a loop of Riemannian metrics on $T^*M$.  One can then consider the gradient of $\ca_H$ with respect to the metric, $\langle \cdot , \cdot \rangle$, written as
$$
\nabla_J\ca_H(z) = -J(z,t)(\frac{dz}{dt}-X_H(t,z)).
$$
The (negative) gradient flow equation on a smooth curve $u : \br \to L(T^*M)$,
$$
\frac{du}{ds} + \nabla_J\ca_H(u(s))
$$
can be rewritten as a perturbed Cauchy-Riemann PDE, if we view $u$ as a smooth map
$\br/\bz \times \br \to T^*M$,  with coordinates, $t \in  \br/\bz, \, s \in \br$,
\begin{equation}\label{cauchy}
\p_s u -J(t,u(t,s))(\p_t u - X_H(t, u(t,s)) = 0.
\end{equation}

Let $a, b \in \cp (H)$. Abbondandolo and Schwarz defined the space of solutions
 
\begin{align}\label{wab}
W(a,b; H,J) = \{u : \br &\to L(T^*M)\, \text{ a solution to} \,  (\ref{cauchy}), \, \text{such that}  \\
\lim_{s\to -\infty}u(s) &= a, \, \text{and} \,  \lim_{s\to +\infty}u(s) = b \}.\notag
\end{align}
 
 As in the case of Morse theory, we then let  $\cm (a,b)$  the ``moduli space" obtained by dividing out by the free $\br$-action,
\begin{equation}\label{flows}
\cm (a,b) = W(a, b; H,J))/\br.
\end{equation}

It was shown in \cite{AS} that with respect to a generic choice of Hamiltonian and almost complex structure, the spaces $W(a, b; H, J)$ and $\cm (a, b)$ whose dimensions are given by $\mu(a) - \mu (b)$ and $\mu (a) - \mu (b) - 1$ respectively, where $\mu$ represents the ``Conley - Zehnder index" of the periodic Hamiltonian orbits $a$ and $b$. Furthermore it was shown in \cite{AS} that in analogy with Morse theory, one can compactify these moduli spaces as
 
$$
\bar \cm (a,b) = \bigcup_{a=a_1 >a_2> \cdots > a_k = b} \cm(a_1, a_2) \times \cdots \times \cm (a_{k-1}, a_k).
$$

The fact that that these compact moduli spaces have canonical framings was shown in \cite{cotangent} using the obstruction to framing described originally in \cite{cjs}.  This was the \sl polarization class \rm defined as follows.

Let  $N$  be an almost complex manifold  whose tangent bundle is classified by a map  $\tau: N \to BU(n)$.   Applying loop spaces, one has a composite map,
\begin{equation}\label{polarization}
\rho: LN \xr{L(\tau)} L(BU(n)) \hk LBU \simeq BU \times U   \to U \to U/O.
\end{equation}
Here the homotopy equivalence $L(BU) \simeq BU \times U$ is well defined up to homotopy, and is given by a trivialization of the   fibration
$$
U \simeq \Omega BU \xr{\iota} L(BU) \xr{ev}BU
$$
where $ev :  LX \to X$ evaluates a loop a $0 \in \br/\bz$.  The trivialization is the composition
$$
U \times BU \xr{\iota \times \sigma} L(BU) \times L(BU) \xr{mult} L(BU).
$$
Here $\sigma : BU \to L(BU) $ is the section of the above fibration given by assigning to a point $x \in BU$ the constant loop at that point, and the  ``multiplication" map in this composition is induced by the infinite loop space structure of $BU$. 

\med
The reason we refer the map $\rho$  as the ``polarization class" of the loop space $LN$, is  because when viewed as an infinite dimensional manifold, the tangent bundle $T(LN)$ is polarized, and its infinite dimensional tangent bundle  has structure group given by the ``restricted general linear group of a Hilbert space", $GL_{res}(H)$ as originally defined in \cite{PS}.  As shown there, $GL_{res}(H)$ has the homotopy type of $\bz \times BO$, and so its classifying space, $BGL_{res}(H) $ has the homotopy type of $B(\bz \times BO) \simeq U/O$ by Bott periodicity.  The classifying map $LN \to U/O$ has the homotopy type of the ``polarization class" $\rho$ defined above.   See \cite{PS} , \cite{cjs}, and \cite{cotangent} for details.  

 When $N = T^*M$ two things were shown in \cite{cotangent}.  First that when viewed as a space of paths, there is a natural map to the based loop space,  $ \iota_{a,b} : \bcm(a,b) \to \Omega L(T^*M)$, well defined up to homotopy, so that the composition
\begin{equation}\label{tangent}
\tau_{a,b} :  \bcm (a, b) \xr{\iota_{a,b}} \Omega L(T^*M) \xr{\Omega \rho} \Omega U/O \simeq \bz \times BO
\end{equation}
  classifies the stable tangent bundle of $\bcm (a,b)$.  Second, it was shown that in this case, i.e when $N = T^*M$, the polarization class $\rho: L(T^*M) \to U/O$ is trivial.  This is essentially because the almost complex structure  (i.e $U(n)$-structure) of the  tangent bundle of $T^*M$ is the complexification of the $n$-dimensional real bundle (i.e $O(n)$-structure) of the tangent bundle of $M$ pulled back to $T^*M$ via the projection map $T^*M \to M$.  By (\ref{tangent}) this leads to a coherent family of framings on the moduli spaces, which in turn lead to a  smooth, compact, framed structure on the flow category $\cc_H$ of the symplectic action functional, as shown in \cite{cotangent}.
 
 Using the methods and results of \cite{AS}, which is to say, comparing the flow category of the perturbed symplectic action functional $\cc_H$ to the Morse flow category of an energy functional on $LM$, the following was shown in \cite{cotangent}.
 
 \med
 \begin{theorem}\label{cotangent}   If $M^n$ is a closed spin manifold of dimension $n$.  For appropriate choices of a Hamiltonian   $H$ and a generic choice of almost complex structure $J$  on the cotangent bundle $T^*M$,   then the Floer homotopy type determined by the smooth, compact flow category $\cc_H$ is given by the suspension spectrum of the free loop space,
 $$
 Z_H(T^*M) \simeq \Sigma^\infty (LM_+).
 $$
 \end{theorem}

 \med
 \noindent
 \bf{Remarks} \rm
 
 1.   This theorem generalized a result of Viterbo \cite{viterbo} stating that the symplectic Floer homology, $SFH_*(T^*M)$ is isomorphic to $H_*(LM)$.
 
 2.  If $M$ is not spin, one needs to use appropriatly twisted coefficients in both Viterbo's theorem and in Theorem \ref{cotangent} above.  This was first observed by Kragh in \cite{kragh}, and was overlooked in all or most of the discussions of the relation between the symplectic Floer theory of the cotangent bundle and homotopy type of the free loop space before Kragh's work, including the author's work in  \cite{cotangent}.

\section{The work of Lipshitz and Sarkar on Khovanov homotopy theory}

A  recent dramatic application of the ideas of Floer homotopy theory appeared in the work of Lipshitz and Sarkar on the homotopy theoretic foundations of Khovanov's homological invariants of knots and links. This work appeared in \cite{lipsark} and \cite{lipsark2}.  Another version of Khovanov homotopy appears in \cite{HKK}.  It was proved to be equivalent to the Lipshitz-Sarkar construction in \cite{LLS}.

The Khovanov homology of a link $L$ is a bigraded abelian group, $Kh^{i,j}(L)$.   It is computed from a chain complex denoted $KhC^{i,j}(L)$ that is defined in terms of a link diagram.  However the Khovanov homology is shown not to depend on the choice of link diagram, and is an invariant of isotopy class of  the link.  This invariant was originally defined by Khovanov in \cite{khov} in which he viewed these homological invariants as a ``categorification" of the Jones polynomial $V(L)$ in the sense that the graded Euler characteristic of this homology theory recovers $V(L)$ via the formula
\begin{align} 
\chi (Kh^{i, j}(L))  &= \sum_{i,j} (-1)^iq^j \, rank \, Kh^{i,j}(L) \notag \\
&= (q + q^{-1})\,  V(L). \notag
\end{align}

\med
The goal of the work of Lipshitz and Sarkar was to associate to a link diagram $L$ a family of spectra $X^j(L)$ whose homotopy types  are invariants of the isotopy class of the link (and in particular do not depend on the particular link diagram used), and so that the Khovanov homology $Kh^{i,j}(L)$ is isomorphic to the reduced singular cohomology $\tilde H^i(X^j(L))$.   Their basic idea is to construct a compact, smooth, framed flow category  from the moduli spaces associated to a link diagram.   Their construction is entirely  combinatorial, and the cells of the spectrum  $X(L)  = \bigvee_j X^j(L)$ correspond to the standard generators of the Khovanov complex $KhC^{*,*} (L)$.  That is,  $X(L)$ \sl realizes \rm the Khovanov chain complex in the sense described above.  $X(L)$ is referred to as the ``Khovanov homotopy type" of the link $L$, and it has had several interesting applications.  

Notice that by virtue of the existence of a Khovanov homotopy type, the Khovanov homology, when reduced modulo a prime, carries an action of the Steenrod algebra $\ca_p$.  In  \cite{lipsark2} the authors show that the Steenrod operation $Sq^2$ acts nontrivially on the Khovanov homology for many knots, and in particular for the torus knot, $T_{3,4}$.  It is also known by work of  Seed \cite{Seed} that there are pairs of links with isomorphic Khovanov's homology, but distinct Khovanov homotopy types.  Also, Rasmussen constructed a slice genus bound, called the $s$-invariant, using Khovanov homology \cite{Ras}. Using the Khovanov homotopy type,  Lipshitz and Sarkar  produced a family of generalizations of the $s$-invariant, and used them to obtain even stronger slice genus bounds.  Stoffregen and Zhang \cite{SZ}  used Khovanov homotopy theory to describe  rank inequalities for Khovanov homology for prime-periodic links in $S^3$. 

\med
We now give a sketch of the construction of compact framed flow category of Lipshitz and Sarkar, which yields the Khovanov homotopy type. 

\med
By a ``link diagram", one means the projection onto $\br^2$ of  an embedded disjoint union of circles in $\br^3$. One keeps track of the resulting ``over" and ``under crossings", and usually one orients the link (i.e puts an arrow in each compoonent).  One can ``resolve" a crossing in two ways.  These are referred to as a $0$-resolution and a $1$ resolution and are described by the following diagram.  

\med
\begin{figure}[ht]
  \centering
  \includegraphics[height=3cm]{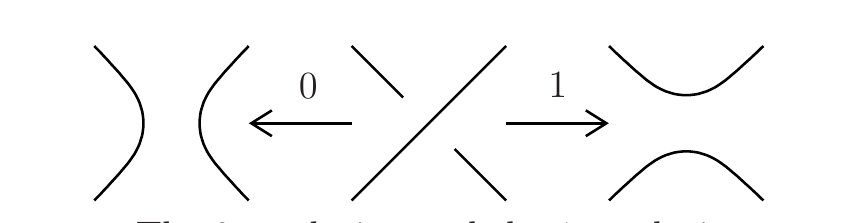}
 \caption{ The $0$-resolution and the $1$-resolution}   
 \label{figone}
\end{figure}

\med
Roughly speaking, the Lipshitz-Sarkar view of the Khovanov chain complex is that it is generated by all possible configurations of resolutions of the crossings of a link diagram.   We recall their definition more carefully.

\med
\begin{definition} A \bf resolution configuration \rm $D$ is a pair $(Z(D), A(D))$, where $Z(D)$ is a set of pairwise disjoint embedded circles in $S^2 = \br^2 \cup \infty$, and $A(D)$ is an ordered collection of arcs embedded in $S^2$ with
$$
A(D) \cap Z(D) = \p A(D).
$$
The number of arcs in the resolution configuration $D$ is called its  \sl index \rm denoted by $ind (D)$.
\end{definition} 

\med
\begin{definition} Given a link diagram $L$ with $n$-crossings, an ordering of the crossings, and a vector $v \in \{0, 1\}^n$, there is an associated resolution configuration $D_L(v)$ obtained by taking the resolution of $L$ corresponding to $v$.  That is, one takes the $0$-resolution of the $i^{th}$ crossing if $v_i = 0$, and the $1$-resolution of the $i^{th}$ crossing if $v_i = 1$.  One then places arcs corresponding to each of the crossings labeled by zero's in $v$.  
\end{definition}

\med
See the following picture of the resolution configuration corresponding to a diagram of the trefoil knot.

\med
\begin{figure}[ht]
  \centering
  \includegraphics[height=5cm]{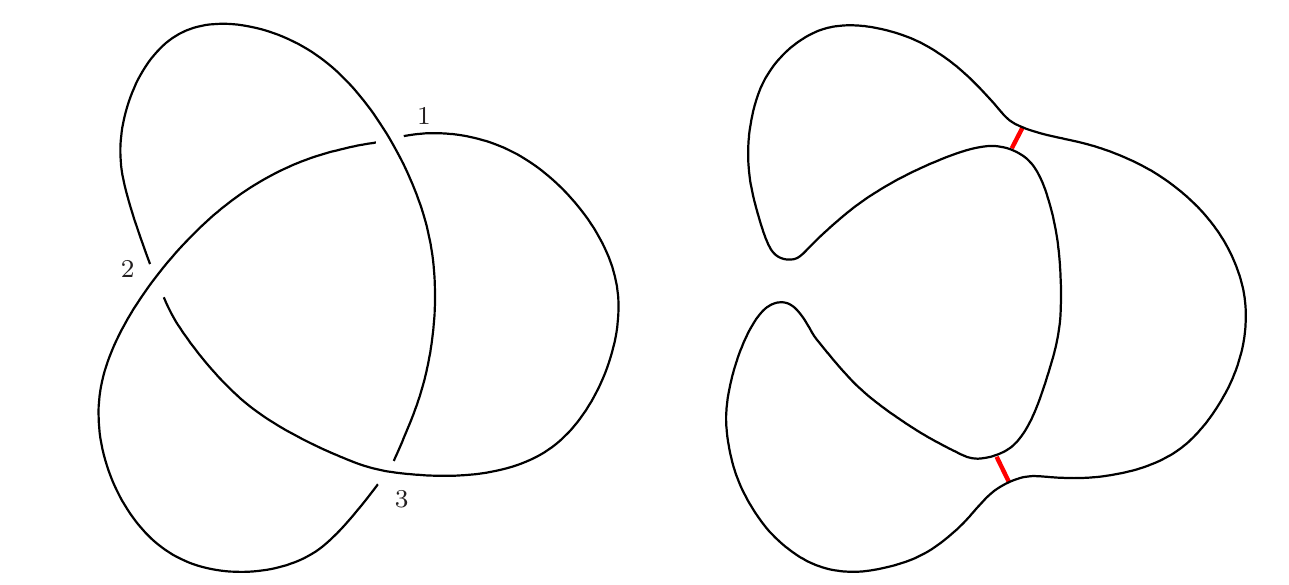}
 \caption{ A knot diagram $K$ for the trefoil with ordered crossings, and  the resolution configuration $D_K((0,1,0))$}   
 \label{figtwo}
\end{figure}

\med
The following terminology is also useful.

\med
\begin{definition}  1.   The \bf core \rm $c(D)$ of a resolution configuration is the the resolution configuration obtained from $D$ by deleting all the circles in $Z(D)$ that are disjoint from all arcs in $A(D0$.

2.  A resolution configuration is \bf basic \rm  if $D = c(D)$, ie every circle in $Z(D)$ intersects an arc in $A(D)$. 
\end{definition}

\med
One can also do a \bf surgery \rm along a subset $A \subset A(D)$. The resulting resolution configuration is denoted $s_{A}(D)$.    The surgery procedure is best illustrated by the following picture.

\begin{figure}[ht]
  \centering
  \includegraphics[height=4cm]{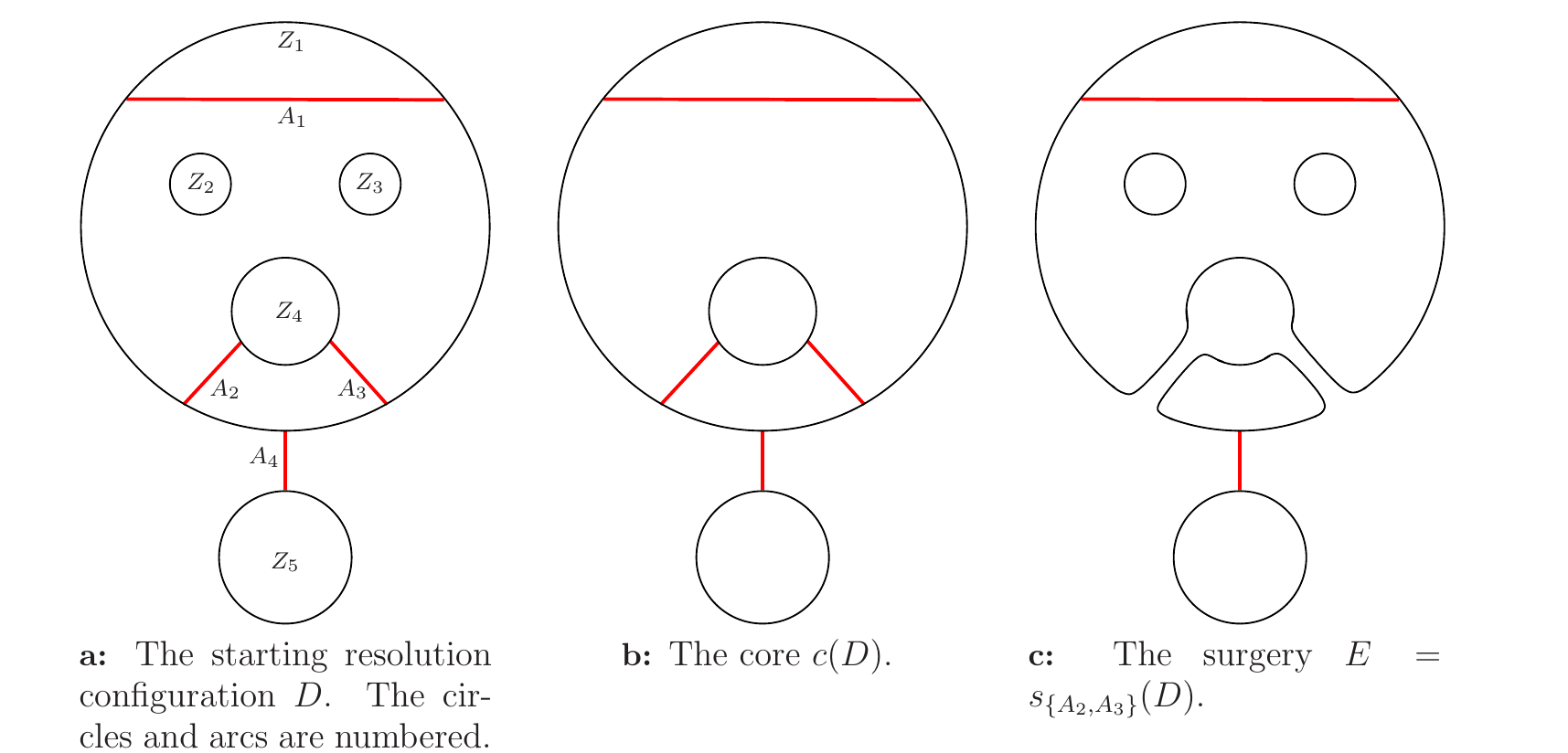}
 %\caption{ A knot diagram $K$ for the trefoil with ordered crossings, and  the resolution configuration $D_K((0,1,0))$}   
 \label{figthree}
\end{figure}

\bg
A \bf labeling \rm $x = \{x_1, x_2\} $ of a resolution configuration $D $  is a labeling of each circle in  $Z:(D)$ by either $x_1$ or $x_2$.    Labeled resolutions configurations have a partial ordering defined to be the transitive closure of the following relations.

We say that $(E,y) < (D,x)$ if
\begin{enumerate}
\item the labelings agree on $D \cap E$
\item $D$ is obtained from $E$ by surgering along a single arc of $A(E)$,  In particular, either:  

(a)  $Z(E \backslash D)$ contains exactly one circle, say $Z_i$ and $Z(D \backslash  E)$ contains exactly two circles, say $Z_j$ and $Z_k$, or

(b) $ Z(E \backslash D)$ contains exactly two circles, say $Z_i$ and $Z_j$, and $Z(D \backslash E)$ contains exactly one circle, say $Z_k$.  
\item
In case (2a), either $ y(Z_i) = x(Z_j) = x(Z_k) = x_- $ or $y(Z_i) = x_+ $ and $\{x(Z_j), x(Z_k)\} = \{x_+, x_-\}$.  

In case (2b), either  $ y(Z_i) = y(Z_j) = x(Z_k) = x_+$ or $y(Z_i) = x_+ $ or $\{y(Z_i), y(Z_j)\} = \{x_+, x_-\}$ and $x(Z_k) = x_+.$
\end{enumerate}

\med
One can now define the \sl Khovanov chain complex \rm $KhC(L)$ as follows. 

\med
\begin{definition} Given an oriented link diagram $L$ with $n$ crossings and an ordering of the crossings in $L$,  $KhC(L)$ is defined to be the free abelian group generated by labeled resolution configurations of the form $(D_L(u), x)$ for $u \in \{0, 1 \}^n$. $ KhC(L)$ carries two gradings, a \sl homological grading \rm $gr_h$ and a \sl quantum grading \rm $gr_q$, defined as follows:
\begin{align} 
gr_h((D_L(u),x)) &= -n_- +|u| \notag \\
gr_q((D_L(u),x)) &= n_+ -2n_- +|u| + \#\{Z \in Z(D_L(u)) \, : \, x(Z) = x_+ \} \notag \\
&- \#\{Z \in Z(D_L(u)) \, : \, x(Z) = x_- \} \notag 
\end{align}
Here $n_+$ denotes the number of positive crossings in $L$, and $n_- $ denotes the number of negative crossings.

The differential preserves the quantum grading, increases the homological grading by $1$, and is defined as
$$
\delta (D_L(v), y) = \sum (-1)^{s_0(\cc_{u,v})} (D_L(u), x)
$$
where the sum is taken over all labeled resolution configurations $(D_L(u), x)$ with $|u| = |v|+1$ and $(D_L(v), y) < (D_L(u), x)$.   The sign $s_0(\cc_{u,v}) \in \bz/2$ is defined as follows:  If $u = (\eps_1, \cdots , \eps_{i-1}, 1, \eps_{i+1}, \cdots , \eps_n)$ and $v = (\eps_1, \cdots , \eps_{i-1}, 0, \eps_{i+1}, \cdots , \eps_n)$, then $s_0(\cc_{u,v}) = \eps_1 + \cdots + \eps_{i-1}.$
\end{definition}

\med
The homology of this chain complex is the \sl Khovanov homology \rm $Kh^{*,*}(L)$.  To define the \sl Khovanov homotopy type \rm of the link $L$, Lipshitz and Sarkar define higher dimensional moduli spaces which have the structure of framed manifolds with corners so that they in turn define a compact, smooth, framed flow category, which by the theory described above, defines the associated (stable) homotopy type.  

These moduli spaces are defined as a certain covering of the moduli spaces occurring in a ``framed flow category of a cube". 
We now sketch these constructions, following Lipshitz and Sarkar \cite{lipsark}.

\med
 Let  $f_1 : \br \to \br$ be a Morse function with one index zero critical point and one index 1 critical point.  For concreteness one can use the function
 $$
 f_1(x) = 3x^2 - 2x^3.
 $$
 Define $f_n : \br^n \to \br$ by
 $$
 f_n(x_1, \cdots , x_n) = f_1(x_1) + \cdots + f_1(x_n).
 $$
 $f_n$ is a Morse function, and we let $\cc (n)$ denote its flow category.  It is a straightforward exercise to see that the geometric realization $|\cc (n)|$ is the $n$-cube $[ 0, 1]^n$.  The vertices of this cube $u \in \{0, 1 \}^n$ correspond to the critical points of $f_n$ and  have a grading which corresponds to the Morse index:
 $gr(u) = |u| = \sum_i u_i$.   They also have a partial ordering coming from the   ordering of $\{0, 1\}$.   We say that $v \leq_i u$ if $v \leq u$ and $gr(u) - gr(v) = i$.  That is $i$ is the relative index of $u$ and $v$.  Let $\bar 1 = (1, 1, \cdots , 1)$ and $\bar 0 = (0, \cdots , 0)$. The following is not difficult, and is verified in \cite{lipsark}.
 
 \med
 \begin{lemma}  The compactified moduli space of piecewise flows $\bar\cm_{\cc (n)} (\bar 1, \bar 0)$ is 
 \begin{itemize}
 \item a single point if $n = 1$, 
 \item a closed interval if $n=2$, 
 \item a closed hexagonal disk if $n = 3$, and is
 \item homeomorphic to a closed disk $D^{n-1}$ for general $n$.
 \end{itemize}
 
 Furthermore, given any  $v < u$, then $\bar \cm_{\cc (n)}(u, v) \cong \bar \cm_{\cc(gr(u)-gr(v))}(\bar 1, \bar 0).$
  \end{lemma}

 \med
 \begin{figure}[ht]
  \centering
  \includegraphics[height=4cm]{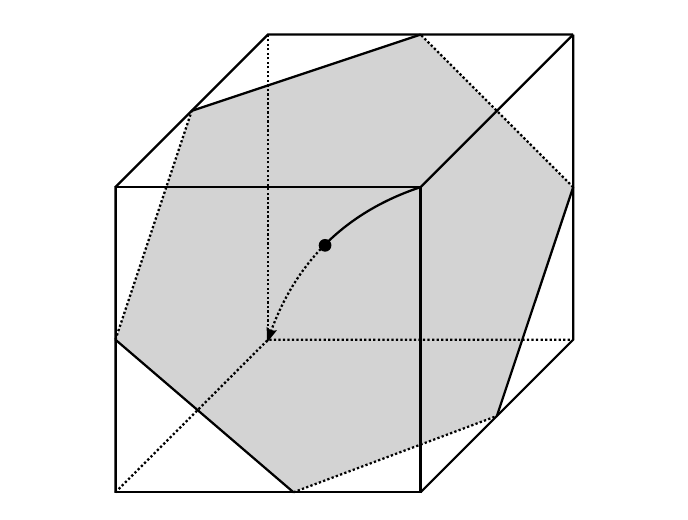}
 \caption{ The cube moduli space  $\bar\cm_{\cc (3)} (\bar 1, \bar 0)$}   
 \label{figfour}
\end{figure}

\med
Lipshitz and Sarkar then proceed to define  moduli spaces of ``decorated resolution configurations" that cover the  moduli spaces occurring in a framed flow category of a cube. 

\med
\begin{definition}  A \bf decorated \rm resolution configuration is a triple $(D, x, y)$ where $D$ is a resolution configuration , $x$ is a labeling of each component of $Z(s(D))$, $y$ is a labeling of each component of $Z(D)$ such that
$$
(D,y) < (s(D), x).
$$
Here $s(D) = s_{A(D)}(D)$ is the maximal surgery on $D$. 
\end{definition}

\med
In \cite{lipsark} Lipshitz and Sarkar proceed to construct moduli spaces $\cm (D,x,y)$ for every decorated resolution configuration $(D, x, y)$. These will be compact, framed manifolds with corners.  Indeed they are $<n-1>$-manifolds where $n$ is the index of $D$.   They also produce  covering maps
$$
\cf : \cm (D, x, y) \to \bar \cm_{\cc (n)} (\bar 1, \bar 0)
$$
which are maps of $<n-1>$-spaces, trivial as a covering maps on each component of $ \cm (D, x, y)$, and are local diffeomorphisms.   The framings of the moduli spaces   $ \cm (D, x, y)$ are then induced from the framings of $ \bar \cm_{\cc (n)} (\bar 1, \bar 0)$.  They produce composition or gluing maps for every labeled resolution configuration $(E,z)$ with 
$$
(D,y) < (E,z) < (s(D), x)
$$
$$
\circ : \cm (D \backslash E, z|, y|) \times \cm (E\backslash s(D), x|, z|) \la  \cm (D, x, y)
$$
that are embeddings into the boundary of $\cm (D, x, y)$.  These moduli spaces are constructed recursively using a clever, but not very difficult argument.  We refer the reader to \cite{lipsark} for details.  

These constructions allow for the definition of a Khovanov flow category for the link $L$, and it is shown to be a compact, smooth, framed flow category as defined  in section 1.2 above.   Using the theory introduced in \cite{cjs} and described in section 1.2,  one obtains a spectrum realizing the Khovanov chain complex.  The  This is the Lipshitz-Sarkar  ``Khovanov homotopy type" of the link $L$. 

\section{Manolescu's equivariant Floer homotopy and the triangulation problem}

In this section the author is relying heavily on the expository article by Manolescu \cite{Man15}.  We refer the reader to this beautiful survey of recent topological  applications of Floer theory.

\subsection{Monopole Floer homology and equivariant Seiberg-Witten stable homotopy}

One example of a dramatic application of Floer's original instanton homology theory, and in particular its topological field theory relationship to the Donaldson invariants of closed $4$-manifolds, was to the study of the group of cobordism classes of homology $3$-spheres.  Define
$$
\Theta^3_H = \{\text{oriented homology $3$-spheres}\}/\sim
$$
where $X_0 \sim X_1$ if and only if there exists a smooth, compact, oriented $4$-manifold $W$ with 
$$
\p W = (-Y_0) \sqcup Y_1
$$ and $H_1(W) = H_2(W) = 0$.  The group operation is represented by connected sum, and the inverse is given by reversing orientation.   The standard unit $3$-sphere $S^3$ is the identity element.   The Rokhlin homomorphism \cite{rok}, \cite{EK} is  the map
$$
\mu  : \Theta^H_3 \to \bz/2
$$
defined by  sending a homology sphere $X$ to $\mu (X) = \sigma (X)/8  \,(mod \, 2)$, where $W$ is any compact spin $4$-manifold with $\p W = X$. It is a theorem that this homomorphism is well-defined. Furthermore, using this homomorphism one knows that the group $\Theta^H_3$ is nontrivial, since, for example, the Poincar\'e sphere $P$ bounds the $E_8$ plumbing which has signature $-8$.  Therefore $\mu (P) = 1$.

This result has been  strengthened using Donaldson theory and Instanton Floer homology.  For example Furuta and Fintushel-Stern proved that $\Theta^H_3$ is infinitely generated \cite{Fur}\cite{FS}.  And using the $SU(2)$ equivariance of Instanton Floer homology, in \cite{Froy}  Froyshov defined a surjective homomorphism 
\begin{equation}\label{froy}
h : \Theta^H_3 \to \bz.
\end{equation}

Monopole Floer homology is similar in nature to Floer's Instanton homology theory, but it is based on the Seiberg-Witten equations rather than the Yang-Mills equations.   More precisely, let $Y$ be a three-manifold equipped with a $Spin^c$ structure $\sigma$.  One considers the configuration space of pairs $(A, \phi)$, where $A$ is a connection on the trivial $U(1)$ bundle over $Y$, and $\phi$ is a spinor.  There is an action of the gauge group of the bundle on this configuration space, and one considers the orbit space of this action.   One can then define the Chern-Simons-Dirac functional on this space by
$$
CSD (A, \phi) = -\frac{1}{8}\int_Y (A^t - A^t_0) \wedge (F_{A^t} + F_{A^t_0})  \, + \frac{1}{2}\int_Y \langle D_A\phi, \phi\rangle \, dvol.
$$
Here $A_0$ is a fixed base connection, and the superscript $t$ denotes the induced connection on the determinant line bundle.  The symbol $F$ denotes the curvature of the connection and the symbol $D$ denotes the covariant derivative.

Monopole Floer homology is the Floer homology associated to this functional.  To make this work precisely involves much analytic, technical work, due in large part to the existence of reducible connections.  Kronheimer and Mrowka dealt with this (and other issues) by considering a blow-up of this configuration space \cite {KM}. They actually defined three versions of Monopole Floer homology for every such pair $(Y, \sigma)$. 

Monopole Floer homology can also be used to give an alternative proof of Froyshov's theorem about the existence of a surjective homomorphism.
$$
\delta : \Theta^H_3 \to \bz
$$
This uses the $S^1$-equivariance of these equations.  Monopole Floer homology has also been used to prove important results in knot theory \cite{KMOS} and in contact geometry \cite{taubes}. 

In \cite{Man03} Manolescu defined a ``Monopole", or ``Seiberg-Witten"  Floer stable homotopy type.  He did not follow the program defined by the author, Jones and Segal in \cite{cjs} as outlined above, primarily because the issue of smoothness in defining a framed, compact, smooth flow category is particularly difficult in this setting.  Instead, he applied Furuta's technique of ``finite dimensional approximations" \cite{Fur2}.   More specifically, the configuration space $X$ of connections and spinors $(A, \phi)$ is a Hilbert space that he approximated by a nested sequence of finite dimensional subspaces $X_\lambda$. He considered the Conley index associated to  the flow induced by $CSD$ on a large ball $B_\lambda \subset X_\lambda$.  Roughly, if $L_\lambda \subset \p B_\lambda$ is that part of the boundary where the flow points in an outward direction, then Manolescu views  the Conley index as the quotient space 
$$
I_\lambda = B_\lambda / L_\lambda.
$$
The homology $H_*(I_\lambda)$ is the Morse-homology of the approximate flow on $B_\lambda$, assuming that the flow satisfies the Morse-Smale transversality condition.   However Manolosecu did not need to assume the Morse-Smale condition in his work.  Namely he simply defined the Seiberg-Witten Floer homology directly as the relative homology of $I_\lambda$, with a degree shift that depends on $\lambda$.  The various $I_\lambda$'s fit together to give a spectrum $SWF(Y, \sigma)$ defined for every rational homology sphere $Y$ with $Spin^c$-structure $\sigma$.  Since the Seiberg-Witten equations have an $S^1$ symmetry, the spectrum $SWF(Y, \sigma)$ carries an $S^1$ action.  This defines Manolescu's ``$S^1$-equivariant Seiberg-Witten Floer stable homotopy type". 

An important case is when $Y$ is a homology sphere.  In this setting there is a unique $Spin^c$ structure $\sigma$ coming from a spin structure.  The conjugation and $S^1$-action together yield an action by the group
$$
Pin (2) = S^1 \oplus S^1j \subset \bc \oplus \bc j = \bh
$$ where $\bh$ are the quarternions.  In \cite{Man13} Manolescu defined the $Pin(2)$-equivariant Floer homology of $Y$ to be the (Borel) equivariant homology of the spectrum,

\begin{equation}\label{Pin}
SWFH_*^{Pin(2)}((Y) = \tilde H^{Pin(2)}_*(SWF(Y)).
\end{equation}

This theory played a crucial role in Manolescu's resolution of the triangulation question as we will describe below.  But  before doing that we point out that by  having the equivariant stable homotopy type, he was able  to define a corresponding $Pin(2)$-equivariant Seiberg-Witten Floer $K$-theory, which he used in \cite{Man13} to prove an analogue of Furuta's   ``10/8"-conjecture for  four-manifolds with boundary.  That is, Furuta proved that if $W$ is a closed, smooth spin four-manifold then 
$$
b_2(W) \geq \frac{10}{8} |sign(W)| + 2
$$
where $b_2$ is the second Betti number and $sign (W)$ is the signature.  (The ``$11/8$-conjecture states  that $b_2 \geq \frac{11}{8} sign(W)$.)   Using $Pin(2)$-equivariant Floer $K$-theory, Manolescu proved that if $W$ is a smooth, spin compact four-manifold with boundary equal to a homology sphere $Y$,  then there is a an invariant $\kappa (Y) \in \bz$, and an analogue of Furuta's inequality,
$$
b_2(W) \geq \frac{10}{8} |sign(W)| + 2 - 2\kappa (Y).
$$

\subsection {The triangulation problem} 

A famous question asked in 1926 by Kneser  \cite{kne} is the following:

\med
\noindent
\bf Question. \rm Does every topological manifold admit a triangulation?

\med
By ``triangulation", one means a homeomorphism to a simplicial complex.  One can also ask the stronger question regarding whether every manifold admits a \sl combinatorial triangulation, \rm which is one in which the links of the simplices are spheres.  Such a triangulation is equivalent to a piecewise linear (PL) structure on the manifold.

In the 1920's Rado proved that every surface admits a combinatorial triangulation, and in the early 1950's, Moise showed that any topological three manifold also admits a combinatorial triangulation.   In the 1930's Cairns and Whitehead showed that \sl  smooth \rm manifolds of any dimension admit combinatorial triangulations.  And in celebrated work in the late 1960's,  Kirby and Siebenmann \cite{KS} showed that there exist topological manifolds without $PL$-structures in every dimension greater than four. Furthermore they showed that in these dimensions, the existence of $PL$-structures is determined by an obstruction class $\Delta (M) \in H^4(M; \bz/2).$   The first counterexample to the simplicial triangulation conjecture was given by Casson \cite{Cas} who showed that Freedman's four dimensional  $E_8$-manifold, which he had proven did not have a $PL$-structure, did not admit a simplicial triangulation.   In dimensions five or greater, a resolution of Kneser's triangulation question
was not achieved until Manolescu's recent work \cite{Man13} using equivariant Seiberg-Witten Floer homotopy theory.

\med
We now give a rough sketch of Manolescu's work on this.

\med
Let $M$ be a closed, oriented $n$-manifold, with $n \geq 5$ that is equipped with a homeomorphism to a simplicial complex $K$ (i.e a simplicial triangulation).  As mentioned above, the Kirby-Siebenmann obstruction to $M$ having a combinatorial ($PL$) triangulation is a class $\Delta (M) \in H^4(M; \bz/2)$.  A related cohomology class is the Sullivan-Cohen-Sato class $c(K)$
\cite{Sul}, \cite{Coh}, \cite{Sato} defined by
$$
c(K) = \sum_{\sigma \in K^{(n-4)} } [link_K(\sigma)] \, \cdot \, \sigma \, \in H_{n-4} (M; \Theta^H_3) \cong H^4(M; \Theta^H_3).
$$
Here the sum is taken over all codimension $4$ simplices in $K$. The link of each such simplex is known to be a homology $3$-sphere.  If this were a combinatorial triangulation the link would be an actual $3$-sphere. 

Consider the short exact sequence given by the Rokhlin homomorphism
\begin{equation}\label{rok}
0 \to ker (\mu) \to \Theta^H_3 \xr{\mu} \bz/2 \to 0.
\end{equation}
This induces a long exact sequence in cohomology
$$
\cdots \to H^4(M; \Theta^H_3) \xr{\mu_*}  H^4(M; \bz/2) \xr{\delta} H^5(M; \ker (\mu_*)) \to \cdots
$$
In this sequence it is known that $\mu_*(c(K)) = \Delta (M) \in H^4(M ;\bz/2).$  One concludes that if $M$ is a manifold that admits a simplicial triangulation, the Kirby-Siebenmann obstruction $\Delta (M)$ is in the image of $\mu_*$ and therefore in the kernel of $\delta$.  An important result was that this necessary condition for admitting a triangulation is also a sufficient condition.
 
 \med
\begin{theorem} \label{GS} (Galewski-Stern \cite{GS}, Matumoto \cite{Mat})    A closed topological manifold $M$ of dimension greater or equal to $5$ is admits a triangulation if and only if $\delta (\Delta (M)) =0$.
\end{theorem}

\med
Now notice that the connecting homomorphism $\delta$ in this long exact sequence would be zero if the short exact sequence (\ref{rok}) were split.  If this were the case then by the Galewski-Stern-Matumoto theorem (\ref{GS}),  every closed topological manifold of dimension $\geq 5$ would be triangulable.  The following theorem states that this is in fact a sufficient condition.

\med

\begin{theorem} \label{GS2} (Galewski-Stern \cite{GS}, Matumoto \cite{Mat})   There exist non-triangulable manifolds of every dimension greater or equal to $5$ if and only if the exact sequence (\ref{rok}) does not split. 
\end{theorem}

\med 
In settling the triangulation problem Manolescu proved the following.

\med
\begin{theorem}\label{triangle} (Manolescu \cite{Man13})  The short exact sequence (\ref{rok}) does not split. 
\end{theorem}

\med
The strategy of his proof was the following.  Suppose $\eta : \bz/2 \to \Theta^H_3$ is a splitting of the above sequence.  The image under $\eta$ of the nonzero element  would represent a homology $3$-sphere $Y$ of order $2$ in $\Theta^H_3$ with nonzero Rokhlin invariant.  Notice that the fact that $Y$ has order $2$ means that $-Y$ ($= Y$ with the opposite orientation) and $Y$ represent the same element in $\Theta^H_3$.  Thus to  show that this cannot happen, Manolescu defined a lift of the Rokhlin invariant to the integers
$$
\beta : \Theta^H_3 \to \bz
$$
Given such a lift and  any element $X \in \Theta^H_3$ of order $2$, 
\begin{align}
 \beta(X) &= \beta(-X) \quad \text{because $X$  has order two} \notag \\
&= -\beta (X) \in \bz   \notag
\end{align}
Thus $ \beta(X) = -\beta (X) \in \bz$ and therefore must be zero. Thus $X$ has zero Rokhlin invariant.

Manolescu's strategy was therefore to construct such a lifting $\beta : \Theta^H_3 \to \bz$ of the Rokhlin invariant $\mu : \Theta^H_3 \to \bz/2$.  His construction was modeled on Froyshov's invariant (\ref{froy}).  But to construct $\beta$ he used $Pin(2)$-equivariant Seiberg-Witten Floer homology $SWFH^{Pin (2)}$ defined using the $Pin(2)$-equivariant Seiberg-Witten Floer homotopy type (spectrum) as described above (\ref{Pin}).  We refer to \cite{Man13} for the details of the construction and resulting dramatic solution of the triangulation problem.

\section{Floer homotopy and Lagrangian immersions: the work of Abouzaid and Kragh}

Another beautiful example of an application of a type of Floer homotopy theory was found by Abouzaid and Kragh \cite{Abou-Kra} in their study of Lagrangian immersions of a Lagrangian manifold $L$ into the cotangent bundle of a smooth, closed manifold, $T^*N$.  Here, $L$ and $N$ must be the same dimension.  The basic question is whether a Lagrangian immersion is Lagrangian isotopic to a Lagrangian embedding.  This is particularly interesting when $L$ = $N$.  In this case the ``Nearby Lagrangian Conjecture" of Arnol'd  states that every closed exact Lagrangian submanifold of $T^*N$ is Hamiltonian  isotopic to the zero section, $\eta : N \hk T^*N$.     In \cite{Abou-Kra} the authors use Floer homotopy theory and classical calculations by Adams and  Quillen \cite{AQ} of the $J$-homomorphism in homotopy theory to  give families of Lagrangian immersions of $S^n$ in $T^*S^n$ that are regularly homotopic to the zero section embedding as smooth immersions, but are \sl not \rm Lagrangian isotopic to any  Lagrangian embedding.

We now describe these constructions and results in a bit more detail.    Let $(M^{2n}, \omega)$ be a symplectic manifold of dimension $2n$.
Recall that a Lagrangian submanifold $L \subset M$ is a smooth submanifold of dimension $n$ such the restriction of the symplectic form $\omega$ to the tangent space of $L$ is trivial. That is, each tangent space of $L$ is an isotropic subspace of the tangent space of $M$.    A Lagrangian embedding $e : L \hk M$ is a smooth embedding whose image is a Lagrangian submanifold.  Similarly, a Lagrangian immersion $\iota : L \to M$ is a smooth immersion so that the image of each tangent space is a Lagrangian subspace of the tangent space of $M$. 

If  $M$ is a cotangent bundle $T^*N$ for some smooth, closed $n$-manifold $N$, then a Lagrangian immersion
$\iota : L \to T^*N$ determines a map $\tau_L : L \to U/O$, which is well-defined up to homotopy.  This map is defined as follows. 

Let $j : N \subset \br^K$ be any smooth embedding with normal bundle $\nu$.  Complexifying, we get
$$
\bc^K \times N \cong (\nu \otimes \bc) \oplus (TN \otimes \bc) \cong (\nu \otimes \bc) \oplus T(T^*N)_{|_N}.
 $$
Here $T(T^*N)$ has a Hermitian structure induced by   its symplectic structure and  the Riemannian structure on N  coming from  the embedding $j$.  So for each $x \in L$,  the tangent space $T_xL$ defines, via the immersion $\iota$, a Lagrangian subspace of $T_{\iota (x)} (T^*N)$. By taking the direct sum with the Lagrangian $\nu  = \nu \otimes \br \subset \nu \otimes \bc$, one obtains a Lagrangian subspace of $\bc^K$.   The Grassmannian of Lagrangian subspaces of $\bc^K$ is homeomorphic to $U(K)/O(K)$.  One therefore has a map
$$
\tau_L : L \to U(K)/O(K) \xr{\subset} U/O.
$$

Now the $h$-principle for Lagrangian immersions states that the set of Lagrangian isotopy classes of immersions $L \to T^*N$  in the homotopy class of a fixed map $\iota_0 : L \to T^*N$ can be identified with the connected components of the space of injective maps of vector bundles
$$TL \to  \iota_0^* (T(T^*N)) $$
which have Lagrangian image.  Let $Sp(T(T^*N))$ be the bundle over $T^*N$ with fiber over  $(x, u) \in T^*N$ the group of linear automorphisms of $T_{(x,u)}(T^*N))$ preserving the symplectic form. This is a principal $Sp(n)$-bundle.  Then the space of all such maps of vector bundles is  either empty or  a principal homogeneous space over the space of sections of $\iota_0^*(Sp(T(T^*N))$.  That is, the space of sections acts freely and transitively.

Abouzaid and Kragh then consider the following special case:

\begin{equation}\label{special}
  L = N,   \,TN \otimes \bc  \,  \text{is a trivial  complex vector bundle, and} \, \iota_0 \, \text{is the inclusion of the zero section.}
   \end{equation}
   
   In this case the the bundle $Sp(T(T^*N)$ is the trivial $Sp(n)$ bundle, and so the corresponding space of sections is the (based) mapping space $Map (N, Sp(n))$, and  since $U(n) \simeq Sp(n)$ it is equivalent to $Map (N, U(n))$.   Since the inclusion of $U(n) \hk U$ is $(2n-1)$- connected,   one concludes the following.
   \begin{lemma}
  The  equivalence classes of Lagrangian immersions in the homotopy class of the zero section of $T^*N$  are classified by homotopy classes of maps from $N$ to $U$.
 \end{lemma}
  
  \med
One therefore has the following homotopy theoretic characterization of the Lagrangian isotopy classes of Lagrangian immersions of the sphere $S^n$  into its cotangent space.

\med
\begin{corollary}\label{sn}
 Isotopy classes of Lagrangian immersions of the sphere $S^n$ in $T^*S^n$  in
the homotopy class of the standard embedding, are classified by  $\pi_n (U)$. 
\end{corollary} 

Of course these homotopy groups are  known  to be the integers $\bz$ when $n$ is odd and zero when $n$ is even.  
  Using a type of Floer homotopy theory, Abouzaid and Kragh proved the following in \cite{Abou-Kra}.
  
  \med
  \begin{theorem}\label{Lagsn}  Whenever $ n$ is congruent to  $1$, $3$, or $5$ modulo $8$, there is a class of Lagrangian immersions of $S^n$  in $T^*S^n$  in the homotopy class of the zero section, which
does not admit an embedded representative.
  \end{theorem}
  
  \med
  We now sketch the Abouzaid-Kragh proof of this theorem, and in particular point out the Floer homotopy  theory  they used.

\med
They first observed that given a Lagrangian immersion $j : N \to T^*N$ in the homotopy class of the zero section, satisfying condition (\ref{special}), then one has   a well defined (up to homotopy)  classifying map $\gamma_j : N \to U$ which lifts the map $\tau_N : N \to U/O$ described above:
\begin{equation}\label{lift}
\tau_N : N \xr{\gamma_j} U \xr{project} U/O.
\end{equation}

Now given any Lagrangian embedding $L \hk M$ as above, Kragh  \cite{kragh} defined a ``Maslov"  (virtual) bundle $\eta_L$ over the  component of the free loop space consisting of contractible loops, $\cl _0L$.  

\med
\noindent \bf Note. \rm We have changed the notation for the free loop space by using a script $\cl$, so as not to get confused by the use of an ``$L$" to denote a Lagrangian. 

\med
The bundle $\eta_L$ is classified by the following map (which by abuse of notation we also call $\eta_L$)
\begin{equation}\label{Maslov}
\eta_L: \cl_0 L \xr{\cl \tau_L} \cl_0 U/O  \xr{\simeq} U/O \times \Omega _0 U/O \xr{project} \Omega_0 U/O \simeq BO.
\end{equation}

Here the equivalence $ \cl_0 U/O \simeq U/O \times \Omega_0 U/O$  comes from considering the evaluation fibration
$$
\Omega U/O \to LU/O \xr{ev} U/O
$$
where $ev$ evaluates a loop at the basepoint.  This fibration 
has a canonical (up to homotopy) trivialization as infinite loop spaces because of the existence of a section as constant loops, and using the infinite loop structure of $U/O$.  The equivalence $\Omega_0U/O \simeq BO$ comes from Bott periodicity.

\med
\noindent \bf Note. \rm  Given a map of any space $f : X \to U/O$, the corresponding virtual bundle over the free loop space
$$
\cl_0 X \to BO
$$
defined as this composition $\cl_0 X \xr{\cl f} \cl_0 U/O \xr{\simeq} U/O \times \Omega _0 U/O \xr{project} \Omega_0 U/O \simeq BO.$ also appeared in the work of Blumberg, the author, and  Schlichtkrull \cite{BCS} in their work on the topological Hochschild homology of Thom spectra.

\med
Let $\bs$ be the sphere spectrum, and let $GL_1(\bs)$ be the group-like monoid of units in $\bs$.  That is, $GL_1(\bs)$  is the colimit of the group-like monoid of based self homotopy equivalences of $S^n$,  $\Omega^n_{\pm 1} S^n$.  These are the degree $\pm 1$ based self maps of $S^n$, which is a group-like monoid under composition.  The classifying space, $BGL_1(\bs)$, classifies (stable) spherical fibrations.  There is a natural map $J : BO \to BGL_1(\bs)$ coming from the inclusion $O(n) \hk \Omega^n_{\pm 1}S^n$ defined by considering the based self equivalence of $ S^n = \br^n \cup \infty $ given by an orthogonal matrix.    We call this map  ``$J$" as it induces the classical $J$ homomorphism on the level of homotopy groups, $ J : \pi_q (O) \to \pi_q(\bs)$.   

Notice that given a $k$-dimensional vector bundle $\zeta \to X$, the associated spherical fibration is the sphere bundle $S(\zeta) \to X$ defined by taking the one-point compactification of each fiber. 

The following is an important result about Lagrangian \sl embeddings \rm  in $\cite{Abou-Kra}.$

\begin{theorem}\label{criterion}  If $j : L \to T^*N$ is an exact Lagrangian embedding then the stable spherical fibration of the Maslov bundle $\eta_L$ is trivial.  That is, the composition
$$
\cl_0 L  \xr{\eta_L} BO \xr{J} BGL_1(\bs)
$$
is null homotopic.
\end{theorem} 

\med
Before sketching how this theorem was proven in \cite{Abou-Kra}, we indicate how Abouzaid and Kragh used this result to detect the Lagrangian immersions of $S^n$ in $ T^*S^n$ yielding Theorem \ref{Lagsn}.  
 
 \med
A Lagrangian immersion $\iota : S^n \to T^*S^n$  is classified by a class $\alpha_\iota \in \pi_n(U)$ by Corollary \ref{sn}.   By Theorem \ref{criterion} if $\iota$  is Lagrangian isotopic to a Lagrangian embedding,  then the composition
 $$
 \cl S^n  = \cl_0 S^n\xr{\cl_0 \alpha_\iota} \cl_0 U  \to \cl_0(U/O) \xr{\simeq} U/O \times \Omega_0 U/O \to \Omega_0U/O  \simeq BO \xr{J} BGL_1(\bs)
 $$
 is null homotopic.  If $\iota$ is such that $\alpha_\iota \in \pi_n(U)$ is a nonzero generator (here $n$ must be odd), then Abouzaid and Kragh precompose this map with the composition
 $$
 S^{n-1} \hk \Omega S^n \to \cl S^n
 $$
 and show that  in the dimensions given in the statement of the theorem, then classical calculations of the $J$-homomorphism in homotopy theory imply  that this composition $S^{n-1} \to BGL_1(\bs)$ is nontrivial.  Therefore by Theorem \ref{criterion}, the Lagrangian immersions of $S^n$ into $T^*S^n$ that these classes represent cannot be Lagrangian isotopic to embeddings, even though as smooth immersions, they are isotopic to the zero section embedding.

 \med
 The proof of Theorem \ref{criterion} is where Floer homotopy theory is used.  This was based on  earlier work by Kragh in \cite{kragh}.  The Floer homotopy theory used was a type of Hamiltonian Floer theory for the cotangent bundle.  They did not directly use the constructions in \cite{cjs} described above, but the spirit of the construction was similar.  More technically they used finite dimensional approximations of the free loop space, not unlike those used by Manolescu \cite{Man03}.  We refer the reader to \cite{kragh}, \cite{Abou-Kra} for details of this construction.   In any case this construction was used to give a spectrum-level version of the Viterbo transfer map, when one has an exact  Lagrangian embedding $j : L \to T^*N$.   This transfer map is given on the spectrum level by a map
 
 $$
 \cl_0 j^! : \cl_0N^{-TN} \to (\cl_0 L)^{-TL \oplus \eta}
 $$
 and similarly a map
 \begin{equation}\label{equivalence}
 \cl_0 j^! : \Sigma^\infty (\cl_0N_+)\to ( \cl_0 L)^{TN -TL \oplus \eta}
 \end{equation}
 An important result in \cite{Abou-Kra} is that $\cl_0j^!$ is an equivalence.  Then they make use of  the  result, essentially due to Atiyah, that given a finite $CW$-complex $X$ with a stable spherical fibration classified by a map $\rho : X \to BGL_1(\bs)$, then the Thom spectrum $X^\rho$ is equivalent to the suspension spectrum $\Sigma^\infty (X_+)$ if and only if $\rho$ is null homotopic.  Applying this to finite dimensional approximations to $\cl_0 L$, they are then able to show that the equivalence (\ref{equivalence}) implies that the Maslov bundle $\eta$, when viewed as a stable spherical fibration is trivial.

 \med
 We end this discussion by remarking on the recasting of the Abouzaid-Kragh results by the author and Klang in \cite{ck}.
 Given a Lagrangian immersion, $j : L \to T^*N$, consider the resulting class 
 $$
 \tau_L : L \to U/O
 $$
 described above.  Taking based loop space, one has a loop map
 $$
 \Omega \tau_L : \Omega_0 L \to \Omega_0 U/O \simeq BO.
 $$
 The resulting Thom spectrum $(\Omega_0 L)^{\Omega \tau_L}$ is a ring spectrum.  Thus one can apply topological Hochshild homology to this ring spectrum,  $THH((\Omega_0 L)^{\Omega \tau_L})$ and one obtains a homotopy theoretic  invariant of the Lagrangian isotopy type of the Lagrangian immersion $j$. This topological Hochshild was computed by Blumberg, the author, and Schlichtkrull in \cite{BCS}.  It was shown that 
 
 $$
 THH((\Omega_0 L)^{\Omega \tau_L}) \simeq (\cl_0L)^{\ell (\tau_L)}
 $$
 where $\ell (\tau_L)$ is a specific stable bundle over the free loop space $\cl L$.  In particular, as was observed in \cite{ck}, in the case when $\tau_L$ factors through a map to $U$, as is the case when $L=N$ and it satisfies the condiion (\ref{special}), then there is an equivalence of stable bundles, $\ell (\tau_L) \cong \eta_L$, where $\eta_L$ is the Maslov bundle  as above.  As a consequence of Theorem \ref{criterion} of Abouzaid and Kragh, one obtains the following:

 \med
 \begin{proposition}  Assume $j : N \to T^*N$ is a Lagrangian immersion in the homotopy class of the zero section, and that the complexification $TN\otimes \bc$ is stably trivial.   Then if, on the level of topological Hochschild homology, we have $
 THH (\Omega_0 N)^{\Omega \tau_N}) ) $ is not equivalent to $  THH(\Sigma^\infty (\Omega N_+)) 
 $,    then $j$ is not Lagrangian isotopic to a Lagrangian embedding.
 \end{proposition}
 
 \med
 Using the results of \cite{BCS}  the authors in \cite{ck}  used this proposition to give a proof of   Abouzaid and Kragh's Theorem \ref{Lagsn}.


\begin{thebibliography}{99}{ 

     \bibitem{AS}A. Abbondandolo and M. Schwarz, \emph{On the Floer homology of cotangent bundles}, Comm. Pure Appl. Math \bf 59, \rm 254-316 (2006) preprint: math.SG/0408280
   
 \bibitem{Abou-Kra} M. Abouzaid and T. Kragh,   \emph{On the immersion classes of nearby Lagrangians},  preprint: arXiv:1305.6810v3 [math.SG]  (2015)
 
 \bibitem{AbS} M. Abouzaid and I. Smith, \emph{ Khovanov homology from Floer cohomology}, arXiv:1504.01230 (2015)
 
  \bibitem{Cas} S. Akbulut and J. McCarthy, \emph{ Casson's invariant for oriented homology $3$-spheres}, Mathematical Notes, vol. 36, Princeton University Press, Princeton, N.J, 1990, an exposition
  
  
 \bibitem{units}M. Ando, A. J. Blumberg, D. J. Gepner, M. J. Hopkins, and C. Rezk.	\emph{Units of ring spectra and Thom spectra}   preprint, arXiv:0810.4535 

\bibitem{BCS} A. Blumberg, R.L Cohen, and C. Schlichtkrull \emph{Topological Hochschild homology of Thom spectra and the free loop space}, Geometry \& Topology \bf vol 14 no.2 \rm{2010}, 1165-1242.


 \bibitem{chassullivan}  M. Chas   and D. Sullivan, \emph{String Topology}.
  preprint: math.GT/9911159.
  
  \bibitem {Coh} M. Cohen,   \emph{Homeomorphisms between homotopy manifolds and their resolutions}, Invent. Math. \bf 10 \rm(1970), 239-250.
  
  
  \bibitem{cotangent}R. L. Cohen, \emph{The Floer homotopy type of the cotangent bundle}, Pure Appl. Math. Q. 6 (2010),
no. 2, Special Issue: In honor of Michael Atiyah and Isadore Singer, 391Ð438. MR 2761853

\bibitem{floeroslo} R. L. Cohen,  \emph{Floer homotopy theory, realizing chain complexes by module spectra, and manifolds with corners,}   in Proc. of Fourth Abel Symposium, Oslo, 2007,  ed. N. Baas, E.M. Friedlander, B. Jahren, P. Ostvaer,  Springer Verlag (2009),  39-59.
preprint: arXiv:0802.2752


\bibitem{cjs1} R.L. Cohen, J.D.S Jones, and G.B. Segal, \emph{Morse theory and classifying spaces}, preprint: http://math.stanford.edu/~ralph/morse.ps (1995)

\bibitem{cjs} R. L. Cohen, J. D. S. Jones, and G. B. Segal, \emph{FloerÕs infinite-dimensional Morse theory and homotopy theory}, The Floer memorial volume, Progr. Math., vol. 133, Birkh¬auser, Basel, 1995, pp. 297Ð
325. MR 1362832 (96i:55012)


\bibitem{ck}R. L. Cohen and I. Klang, \emph{ 
Twisted Calabi-Yau ring spectra, string topology, and gauge symmetry}, preprint  arXiv:1802.08930 (2018)

 

%  \bibitem{cohenjones} R.L. Cohen and J.D.S. Jones, \emph{ A homotopy theoretic realization 
%of string topology},  
%  Math. Annalen, \bf  vol. 324,  \rm 773-798 (2002).  preprint: math.GT/0107187 


\bibitem{EK} J. Eels, Jr. and N. Kuiper,  \emph{An invariant for certain smooth manifolds} Ann. Mat. Pura Appl. (4) \bf{60} \rm (1962), 93-110.
 
 
 \bibitem{FS} R. Fintushel and R. Stern, \emph{Instanton homology of Seifert fibred homology three spheres} Proc. London Math. Soc. (3) \bf{61} \rm (1990), no. 1, 109-137.
 
\bibitem{floerLag} A. Floer,  \emph{Morse theory for Lagrangian intersections,} J. Differential Geometry \bf 28 \rm (1988),  513-547.

\bibitem{floerI} A. Floer, \emph{An instanton invariant for 3-manifolds,} Communications in Mathematical Physics, \bf 118 \rm (1988), 215-240.

\bibitem{floerFix} A. Floer, \emph{Symplectic fixed points and holomorphic spheres,} Communications in Mathematical Physics, \bf   120 \rm (4) (1989), 575Ð611. 


\bibitem{franks}  J. M. Franks, \emph{Morse-Smale flows and homotopy theory}, Topology \bf 18 \rm (1979), 119-215.

\bibitem{Froy} K. A. Froyshov, \emph{Equivariant aspects of Yang-Mills Floer theory}, Topology \bf{41} \rm (2002), no. 3, 525-552.

\bibitem{Fur} M. Furuta, \emph{Homology cobordism group of homology $3$-spheres,} Invent. Math.  \bf 100 \rm (1990), no. 2, 339-355.

\bibitem{Fur2} M. Furuta, \emph{Monopole equation and the $\frac{11}{8}$-conjecture}, Math. Res. Lett. \bf 8 \rm (2001), no. 3, 279-291.

\bibitem{GS} D. Galewski and R. Stern, \emph{Classification of simplicial triangulations of topological manifolds}, Ann. Math. (2) \bf 111 \rm (1980), no.1, 1-34.

\bibitem{genauer} J. Genauer,  \emph{Cobordism categories of manifolds with corners}, Stanford University PhD thesis,  (2009), preprint https://arxiv.org/abs/0810.0581

\bibitem{HKK} P. Hu, D. Kriz, and I. Kriz, \emph{Field theories, stable homotopy theory, and Khovanov homology}, arXiv:1203.4773 (2012)

\bibitem{khov} M. Khovanov, \emph{A categorification of the Jones polynomial},  Duke Math. J. \bf 101 \rm no. 3, (2000), 
359Ð426.

\bibitem{KS} R. Kirby and L. Siebenmann, \emph{Foundational Essays on topological manifolds, smoothings, and triangulations}, Princeton University Press, Princeton, N.J., 1977, with notes by J. Milnor and M. Atiyah, Annals of Math. Studies, No. 88.

\bibitem{kne} H. Kneser, \emph{Die Topologie der Mannigfaltigkeiten}, Jahresbericht der Deut. Math. Verein. \bf 34 \rm (1926), 1-13.

\bibitem{kragh}  T. Kragh, \emph{Parameterized ring-spectra and the nearby lagrangian conjecture}, Geom \& Top. \bf 17 \rm (2013), 639-731 (Appendix by M. Abouzaid)


\bibitem{KM} P. Kronheimer and T. Mrowka, \bf Monopoles and three-manifolds, \rm New Mathematical Monographs, vol. 10, Cambridge University Press, Cambridge, 2007. 

\bibitem{KMOS} P. Kronheimer, T. Mrowka, P. Oszvath, and Z. Szabo, \emph{Monopoles and lens space surgeries}, Ann. of Math. (2) \bf 165 \rm (2007), no. 2, 457-546.

\bibitem{laures}G. Laures, \emph{On cobordism of manifolds with corners}, Trans. of AMS \bf 352 \rm no. 12,  (2000), 5667-5688.

\bibitem{LLS}T. Lawson, R. Lipshitz, and S. Sarkar, \emph{Khovanov homotopy type, Burnside category, and products}, arXiv: 1505.00213v1 (2015).

\bibitem{lipsark} R. Lipshitz and S. Sarkar, \emph{A Khovanov stable homotopy type}, J. Amer. Math. Soc. 27 (2014), no. 4, 983-1042.

\bibitem{lipsark2} R. Lipshitz and S. Sarkar, \emph{A Steenrod square on Khovanov homology}, J. Topology, 7 (2014), no. 3, 817-848. arXiv:1204.5776.

\bibitem{Man03} C. Manolescu, \emph{Seiberg-Witten-Floer stable homotopy type of three-manifolds with b1 = 0},
Geom. Topol. 7 (2003), 889Ð932 (electronic). MR 2026550 (2005b:57060)

\bibitem{Man07} C. Manolescu, \emph{A gluing theorem for the relative Bauer-Furuta invariants}, J. Differential Geom. 76 (2007),
no. 1, 117Ð153. MR 2312050 (2008e:57033)


\bibitem{Man13} C. Manolescu, \emph{Pin (2)-equivariant Seiberg-Witten Floer homology and the triangulation conjec-
ture}, arXiv preprint arXiv:1303.2354 (2013).

\bibitem{Man15} C. Manolescu, \emph{Floer theory and its topological applications}, preprint arXiv:1508.00495 (2015).

\bibitem{Man16} C. Manolescu, \emph{The Conley index, gauge theory, and triangulations}, arXiv preprint arXiv:1308.6366v4, (2016).

\bibitem{Mat} T. Matumoto, \emph{Triangulations of manifolds}, Algebraic and Geometric Topology (Proc. Symp. Pure Math., Stanford Univ., 1976), Part 2, Proc. Symp. Pure Math., XXXII, AMS, 1978, pp.3 - 6.

\bibitem {OS} P.  Ozsvath and Z.  Szabo  \emph{  Holomorphic disks and topological invariants for closed three-manifolds}  Annals of Mathematis \bf 159  \rm(3) (2004)  1027Ð1158
  \bibitem {OS2} P.  Ozsv‡th and Z.  Szabo  \emph{  Holomorphic disks and knot invariants}, arXiv:math/0209056, (2003).
  
  \bibitem{PS} 
 A. Pressley and G. Segal,
\textbf{ Loop Groups},
  Oxford Math. Monographs, Clarendon Press
(1986).

\bibitem{Ras} J. Rasmussen,  \emph{ Khovanov homology and the slice genus},  Invent. Math. 182 (2010), no. 2, 419Ð447.
 
 \bibitem{AQ} D. Ravenel, \bf Complex Cobordism and Stable Homotopy Groups of Spheres, \rm  Pure and Applied Mathematics, vol. 121, Academic Press, Inc., Orlando, FL, 1986
 
 
 \bibitem{rok} V. A. Rokhlin, \emph{New results in the theory of four-dimensional manifolds}, Doklady Akad. Nauk SSSR (N.S) {\bf 84} (1952), 221-224.
  
  \bibitem{qin} L. Qin,    \emph{On the Associativity of Gluing},            preprint:  http://front.math.ucdavis.edu/1107.5527
   
     \bibitem{Ras} J.  Rasmussen, \emph{Floer homology and knot complements}, 	arXiv:math/0306378, (2003)

\bibitem{Sato} H. Sato, \emph{Constructing manifolds by homotopy equivalences, I. An obstruction to constructing PL manifolds from homology manifolds}, Ann. Inst. Fourier (Grenoble) \bf 22 \rm (1972), no. 1, 271-286.

\bibitem{Seed} C. Seed, \emph{Computations of the Lipshitz-Sarkar Steenrod square on Khovanov homology},
arXiv:1210.1882.

\bibitem{SZ} M. Stoffregen and M. Zhang, \emph{Localization in Khovanov homology}, arXiv:1810.04769, (2018).

\bibitem{SS} P. Seidel and I. Smith, \emph{A link invariant from the symplectic geometry of nilpotent slices},  Duke Math. J.
\bf 134  \rm (2006),  453-514.

\bibitem{Sul} D. Sullivan, \emph{Triangulating and smoothing homotopy equivalences and homeomorphisms}, Geometric Topology Seminar notes, The Hauptvermutung book, $K$-Monogr. Math., vol.1, Kluwer Acad. Publ.,Dordrecht, 1996, pp. 69-103.

\bibitem{taubes} C. Taubes \emph{The Seiberg-Witten equations and the Weinstein conjecture}, Geom. \& Topol. \bf 11 \rm (2007), 2117 - 2202.
         \bibitem{viterbo}C. Viterbo, \emph{Functors and computations in Floer homology with applications, Part II}, preprint, (1996). 
         
                  }\end{thebibliography}
\end{document}